\documentclass{gtart}     % Basic GT style
\usepackage{amsmath,amssymb,amscd,epsf}
\input gtmonout
\volumenumber{2}
\volumeyear{1999}
\volumename{Proceedings of the Kirbyfest}
\pagenumbers{259}{290}
\papernumber{14}
\received{28 December 1998}\revised{30 May 1999}
\published{20 November 1999}

\newtheorem{theorem}{Theorem}[section]
\newtheorem{lemma}[theorem]{Lemma}
\newtheorem{corollary}[theorem]{Corollary}
\newtheorem{proposition}[theorem]{Proposition}

\theoremstyle{definition}
\newtheorem{definition}[theorem]{Definition}
\newtheorem{example}[theorem]{Example}

\theoremstyle{remark}

\numberwithin{equation}{section}
\let\Bbb\mathbb

\begin{document}
\title{Classification of unknotting tunnels\\for two bridge knots}        
\authors{Tsuyoshi Kobayashi}                  
\address{Department of Mathematics, Nara
Women's University\\Kita-Uoya Nishimachi, Nara 630, JAPAN}                  
\email{tsuyoshi@cc.nara-wu.ac.jp}                     
\begin{abstract} 
In this paper, we show that any unknotting tunnel 
for a two bridge knot is isotopic to either one 
of known ones. 
This together with Morimoto--Sakuma's result gives the 
complete classification of unknotting tunnels for 
two bridge knots up to isotopies and homeomorphisms. 
\end{abstract}

\makeshorttitle

\section{Introduction}

Let $K$ be a knot in the 3--sphere $S^3$. 
The {\it exterior} of $K$ is the closure of the 
complement of a regular neighborhood of $K$, 
and is denoted by $E(K)$.  
A {\it tunnel} for $K$ is an embedded arc $\sigma$ 
in $S^3$ such that $\sigma \cap K = \partial \sigma$.  
Then we denote $\sigma \cap E(K)$ by $\hat \sigma$, 
where we regard $\sigma$ as obtained from 
$\hat \sigma$ by a radial extension. 
Let $\sigma _1$, $\sigma _2$ be tunnels for $K$.  
We say that $\sigma _1$ and $\sigma _2$ are 
{\it homeomorphic} if there is a self homeomorphism 
$f$ of $E(K)$ such that $f(\hat \sigma_1) = \hat \sigma_2$.  
We say that 
$\sigma_1$ and $\sigma_2$ are {\it isotopic} 
if $\hat \sigma_1$ is ambient isotopic to $\hat \sigma_2$ 
in $E(K)$.  

We say that a tunnel $\sigma$ for $K$ is {\it unknotting} 
if $S^3-$Int~$N(K \cup \sigma ,S^3)$ 
is a genus two handlebody.  
We note that the unknotting tunnels for $K$ is 
essentially the genus 2 Heegaard splittings of 
$E(K)$; 
if $\sigma$ is an unknotting tunnel, 
then 
we can obtain a genus 2 Heegaard 
splitting $(C_1, C_2)$, 
where 
$C_1$ is a regular neighborhood of 
$\partial E(K) \cup \hat \sigma$ in $E(K)$, 
and 
$C_2 = c\ell(E(K)-C_1)$, 
and every genus 2 Heegaard splitting of 
$E(K)$ is obtained in this manner. 
Moreover, such Heegaard splittings 
are isotopic (homeomorphic resp.) 
if and only if 
the corresponding unknotting tunnels are 
isotopic (homeomorphic resp.). 
We say that a knot $K$ is a 2--bridge knot 
if $K$ admits a (genus zero) 2--bridge position, 
that is, 
there exists a genus zero Heegaard splitting 
$B_1 \cup _P B_2$ of 
$S^3$ such that $K \cap B_i$ is a system 
of 2--string trivial arcs 
in $B_i$ $(i=1,2)$.
It is known that each 2--bridge knot admits six 
unknotting tunnels as depicted in Figure~1.1 or 
Figure~3.1 (see \cite{T}, or \cite{K}).

\begin{figure}[ht!]\small
\begin{center}
\leavevmode
\epsfxsize=90mm
\epsfbox{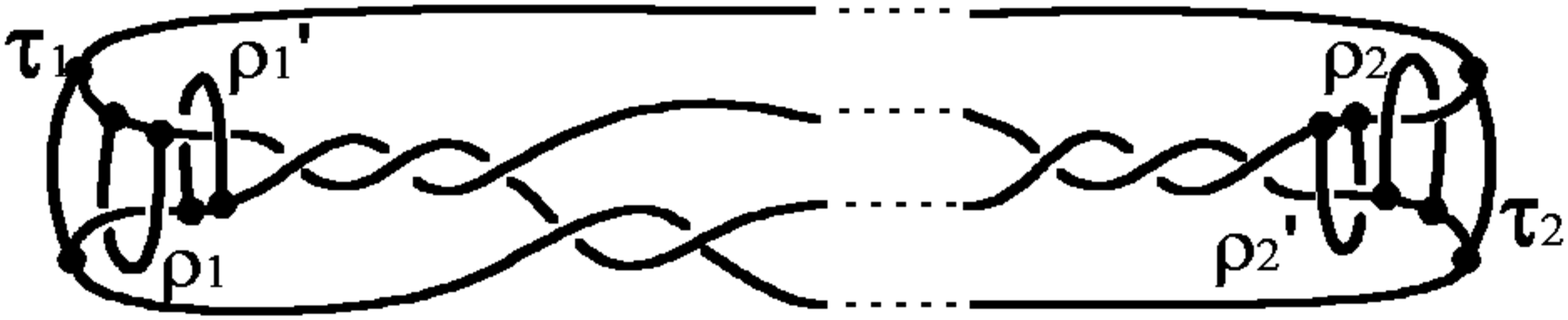}
\end{center}

\begin{center}
Figure 1.1
\end{center}
\end{figure}

Then the purpose of this paper is to prove: 

\begin{theorem}\label{main}
Every unknotting tunnel for a non-trivial 2--bridge knots
is isotopic to one of the above six unknotting tunnels. 
\end{theorem}

We note that the isotopy, and 
homeomorphism classes of the above tunnels are 
completely classified by Morimoto--Sakuma \cite{M-S} 
and Y.Uchida \cite{U}, and 
that it is known that the unknotting tunnels for 
a trivial knot are mutually isotopic (see, for example 
\cite{S-T}). Hence these results together with the above 
theorem give the complete classification of 
isotopy, and homeomorphism classes of unknotting 
tunnels for two--bridge knots.

\section{Preliminaries}

Throughout this paper, 
we work in the differentiable category. 
For a submanifold $H$ of a manifold $K$, 
$N(H;K)$ denotes a regular neighborhood of $H$ in $K$.
Let $N$ be a manifold embedded in a manifold $M$
with dim$N=$dim$M$.
Then $\text{Fr}_M N$ denotes the frontier of $N$ in $M$.
For the definitions of standard terms in 
3--dimensional topology, 
we refer to \cite{He}.

Let $M$ be a compact 3--manifold, 
$\gamma$ a union of mutually disjoint arcs 
or simple closed curves properly 
embedded in $M$, 
$F$ a surface embedded in $M$, which is in general position 
with respect to $\gamma$, 
and $\ell (\subset F)$ a simple closed curve with 
$\ell \cap \gamma = \emptyset$. 

\begin{definition}
A surface $D$ in $M$ is a {\it $\gamma$--disk}, 
if $D$ is a disk intersecting $\gamma$ in 
at most one transverse point. 
\end{definition}

\begin{definition}
We say that $\ell$ is $\gamma$--{\it inessential} 
if $\ell$ bounds a $\gamma$--disk in $F$. 
We say that $\ell$ is $\gamma$--{\it essential} if 
it is not $\gamma$--inessential.
\end{definition}

\begin{definition}
We say that a disk $D$ is a $\gamma$--{\it compressing disk} 
for $F$ if; 
$D$ is a $\gamma$--disk, and 
$D \cap F = \partial D$, and 
$\partial D$ is a $\gamma$--essential simple closed curve in $F$. 
The surface $F$ is $\gamma$--{\it compressible} 
if it admits a $\gamma$--compressing disk, and 
it is $\gamma$--{\it incompressible} if it is not 
$\gamma$--compressible. 
\end{definition}

\begin{definition}
Let $F_1$, $F_2$ be 
surfaces embedded in $M$ 
such that $\partial F_1 = \partial F_2$, 
or 
$\partial F_1 \cap \partial F_2 = \emptyset$. 
We say that $F_1$ and $F_2$ are $\gamma$--{\it parallel}, 
if there is a submanifold $N$ in $M$ such that 
$(N, N \cap \gamma )$ is homeomorphic to 
$(F_1 \times I, {\cal P} \times I)$ 
%as triple, where ${\cal P}$ is a union of points in
as a pair, where ${\cal P}$ is a nion of points in
$\text{Int}(F_1)$, and 
%$F_1$ ($F_2$ resp.) is contained in the 
$F_1$ ($F_2$ resp.) is the 
closure of the component of 
$\partial (F_1 \times I) - (\partial F_1 \times \{ 1/2 \})$ 
containing $F_1 \times \{ 0 \}$ ($F_1 \times \{ 1 \}$ resp.) 
if $\partial F_1 = \partial F_2$, or 
$F_1$ ($F_2$ resp.) is the surface corresponding to 
$F_1 \times \{ 0 \}$ 
($F_1 \times \{ 1 \}$ resp.) if 
$\partial F_1 \cap \partial F_2 = \emptyset$.

The submanifold $N$ is called a {\it $\gamma$--parallelism} 
between $F_1$ and $F_2$. 

We say that $F$ is $\gamma$--{\it boundary parallel} 
if there is a subsurface $F'$ in $\partial M$  such that 
$F$ and $F'$ are $\gamma$--parallel. 
\end{definition}

\begin{definition}
We say that $F$ is {\it $\gamma$--essential} 
if $F$ is $\gamma$--incompressible, and 
not $\gamma$--boundary parallel. 
\end{definition}

Let $a$ be an arc properly embedded in $F$ 
with $a \cap \gamma = \emptyset$.

\begin{definition}
We say that $a$ is {\it $\gamma$--inessential} 
if there is a subarc $b$ of $\partial F$ 
such that $\partial b = \partial a$, 
and $a \cup b$ bounds a disk $D$ in $F$ such that 
$D \cap \gamma = \emptyset$, 
and $a$ is {\it $\gamma$--essential} if it is not 
$\gamma$--inessential. 
\end{definition}

\begin{definition}
We say that $F$ is {\it $\gamma$--boundary compressible} 
if there is a disk $\Delta$ in $M$ such that 
$\Delta \cap F = \partial \Delta \cap F = \alpha$ 
is an $\gamma$--essential arc in $F$, and 
$\Delta \cap \partial M = \partial \Delta \cap \partial M = 
c\ell (\partial \Delta - \alpha)$. 
\end{definition}

\begin{definition}
Let $F_1$, $F_2$ be mutually disjoint surfaces in $M$ 
which are in general position with respect to $\gamma$. 
We say that $F_1$ and $F_2$ are {\it $\gamma$--isotopic} 
if there is an ambient isotopy 
$\phi_t$ $(0 \le t \le 1)$ of $M$ such that; 
$\phi_0 = id_M$; 
$\phi_1(F_1) = F_2$, and; 
$\phi_t(\gamma ) = \gamma$ for each $t$. 
\end{definition}

{\bf Genus $g$ $n$--bridge position}

Let $\Lambda  = \{ \gamma_1, \dots , \gamma_n \}$ be a system 
of mutually disjoint arcs properly embedded in $M$. 

\begin{definition}
We say that $\Lambda$ is a 
{\it system of $n$--string trivial arcs} 
if there exists 
a system of mutually disjoint disks 
$\{ D_1, \dots , D_n \}$ in $M$ 
such that, for each $i$ $(i = 1, \dots , n)$, we have 
(1) $D_i \cap \Lambda = \partial D_i \cap \gamma_i 
= \gamma_i$, and 
(2) $D_i \cap \partial M$ is an arc, say 
$\alpha_i$, such that $\alpha_i = c\ell(\partial D_i - \gamma_i)$. 
\end{definition}

\begin{example}
Let $\beta$ be a system of 2--string trivial arcs 
in a 3--ball $B$. 
The pair $(B, \beta )$ is often refered as {\it 2--string trivial tangle}, 
or a {\it rational tangle}.
\end{example}

Let $K$ be a link in a closed 3--manifold $M$. 
Let $M = A \cup_P B$ be a genus $g$ Heegaard splitting. 
Then the next definition is borrowed from \cite{D}. 

\begin{definition}
We say that $K$ is in a {\it (genus $g$) $n$--bridge position} 
(with respect to the Heegaard splitting $A \cup_P B$) if 
$K \cap A$ ($K \cap B$ resp.) is a system of 
$n$--string trivial arcs 
in $A$ ($B$ resp.). 
\end{definition}

In this paper, we abbreviate a genus 0 $n$--bridge 
position to an $n$--bridge position. 
A knot $K$ is called an {\it $n$--bridge knot} 
if it admits an $n$--bridge position. 
It is known that the $2$--bridge positions of 
a 2--bridge knot $K$ are unique up to 
$K$--isotopy (see \cite{O},\cite{Sch}, or Section 7 of \cite{K-S}). 

\begin{definition}\label{2-bridge K-reducible}
We say that a genus $g$ bridge position of $K$ 
with respect to $A \cup_P B$ is 
{\it weakly $K$--reducible} 
if there exist $K$--compressing disks $D_A$, $D_B$ for $P$ in 
$A$, $B$ respectively 
such that $\partial D_A \cap \partial D_B = \emptyset$. 
The genus $g$ bridge position of $K$ 
with respect to $A \cup_P B$ is 
{\it strongly $K$--irreducible} 
if it is not weakly $K$--reducible.
\end{definition}

\medskip
\noindent
{\bf Remark}\qua
It is known that the 2--bridge positions of a 2--bridge knot 
are strongly $K$--irreducible 
(see Proposition~7.5 of \cite{K-S}). 

\medskip
For a 2--bridge knot $K$ we can obtain four genus one 1--bridge 
positions of $K$ as follows. 

\medskip
Let $A \cup_P B$ be the Heegaard splitting which 
gives the 2--bridge position, and 
$a_1$, $a_2$, $b_1$, $b_2$ the 
closures of the components of $K - P$, 
where $a_1 \cup a_2$ ($b_1 \cup b_2$ resp.) 
is contained in $A$ ($B$ resp.). 
Let 
$T_1 = A \cup N(b_1,B)$, 
$\alpha_1 = a_1 \cup b_1 \cup a_2$, 
$T_2 = c\ell(B -N(b_1,B))$, and 
$\alpha_2 = b_2$. 
Then each $T_i$ is a solid torus and 
it is easy to see that 
$\alpha_i$ is a trivial arc in $T_i$ 
($i = 1,2$). 
Hence, $T_1 \cup T_2$ gives genus one 
1--bridge position of $K$. 
Moreover, by using $a_1$, $a_2$, $b_2$ for $b_1$, 
we can obtain other three genus one 1--bridge 
positions of $K$.

\medskip
Let $K$ be a knot with a genus one 1--bridge 
position with respect to $T_1 \cup T_2$. 
Let 
$\mu_1$, $\mu_2$ be tunnels for $K$ 
embedded in $T_1$, $T_2$ respectively 
as in Figure~2.1. 
It is easy to see that $\mu_1$, $\mu_2$ 
are unknotting tunnels, and we call them the 
{\it unknotting tunnels 
associated to the genus one 1--bridge position}. 
In Section 8 of \cite{K-S}, it is shown that 
every genus one 1--bridge position for a non-trivial 
2--bridge knot is obtained as above. 
Hence, by definition (see also Figure~3.1), it is easy 
to see: 

\begin{proposition}\label{1to2}
Let $\mu_1$, $\mu_2$ be unknotting tunnels 
associated to a 
genus one 1--bridge position of a 2--bridge knot $K$. 
Then one of $\mu_1$, $\mu_2$ is isotopic to $\tau_1$ or $\tau_2$, 
and 
the other is isotopic to either 
$\rho_1$, $\rho_1'$, 
$\rho_2$ or $\rho_2'$ 
\end{proposition}

\begin{figure}[ht!]\small
\begin{center}
\leavevmode
\epsfxsize=70mm
\epsfbox{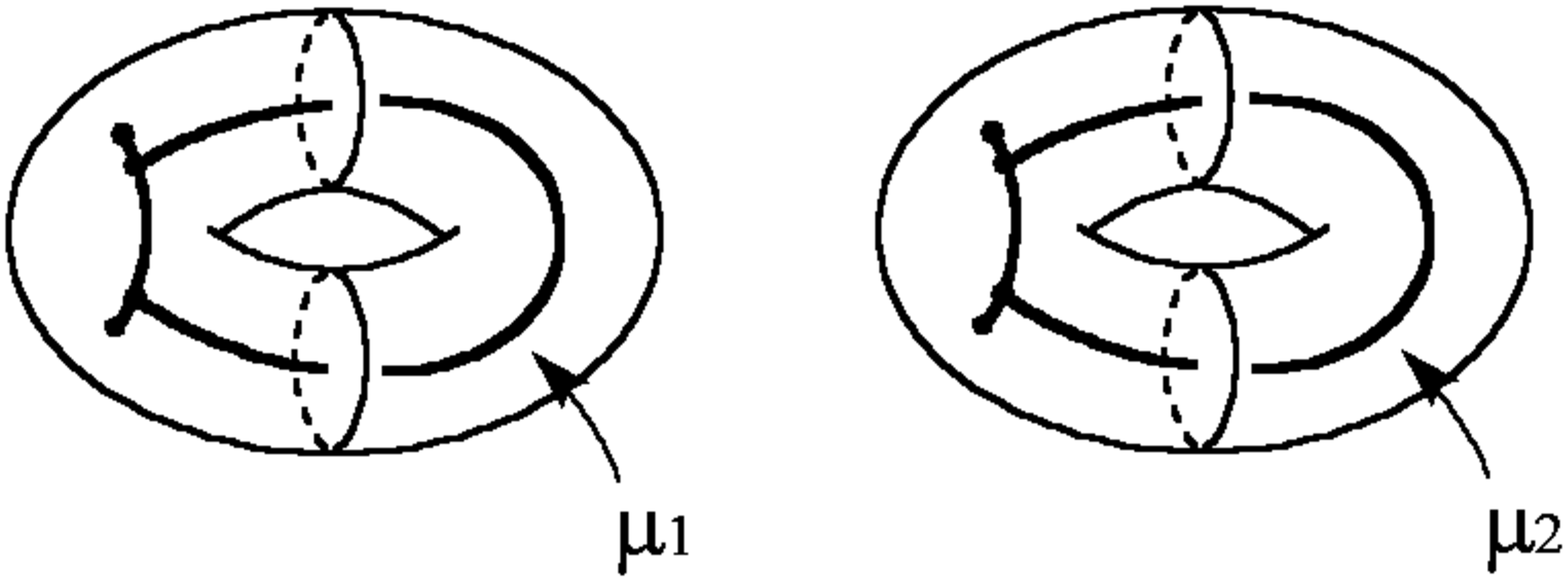}
\end{center}
\begin{center}
Figure~2.1 
\end{center}
\end{figure}

Let $\sigma$ be an unknotting tunnel for $K$. 
Let 
$V_1 = N(K \cup \sigma ; S^3)$, 
$V_2 = c\ell(S^3-V_1)$. 
Note that 
$V_1 \cup_Q V_2$ is a genus two Heegaard splitting of $S^3$. 
%%%%%%

\begin{definition}
We say that the Heegaard splitting 
$V_1 \cup V_2$ is 
{\it weakly $K$--reducible} 
if there exist $K$--compressing disks $D_1$, $D_2$ 
properly embedded in 
$V_1$, $V_2$ respectively 
such that $\partial D_1 \cap \partial D_2 = \emptyset$. 
The splitting is
{\it strongly $K$--irreducible} 
if it is not weakly $K$--reducible.
\end{definition}

\begin{proposition}\label{weakly reducible}
If $(V_1, V_2)$ is weakly $K$--reducible, 
then either $K$ is a trivial knot or 
$K$ admits a genus one 1--bridge position, 
where $\sigma$ is isotopic to one of the unknotting tunnels 
associated to the 1--bridge position. 
\end{proposition}

\begin{proof}
Let $D_1$ ($\subset V_1$), $D_2$ ($\subset V_2$) 
be a pair of $K$--compressing disks 
which gives weak $K$--irreducibility. 

\medskip
\noindent
{\bf Claim 1}\qua 
We may suppose that $D_1$ ($D_2$ resp.) is non-separating 
in $V_1$ ($V_2$ resp.). 

\medskip
\noindent
{\bf Proof of Claim 1}\qua
Suppose that $D_2$ is separating in $V_2$. 
Then $D_2$ cuts $V_2$ into two solid tori, 
say $T_1$, $T_2$. 
By exchanging the suffix, if necessary, 
we may suppose that 
$\partial D_1 \subset \partial T_1$. 
Then take a meridian disk 
$D_2'$ in $T_2$ such that 
$\partial D_2' \subset \partial V_2$. 
We may regard $D_2'$ as a 
(non-separating essential) disk in $V_2$, 
and we have 
$\partial D_1 \cap \partial D_2' = \emptyset$. 
By regarding $D_2'$ as $D_2$, we see that 
we may suppose that $D_2$ is non-separating in $V_2$. 

Suppose that 
$D_1$ is separating in $V_1$. 
Since $K$ does not intersect $D_1$ in one point, 
we have $D_1 \cap K = \emptyset$. 
The disk $D_1$ cuts $V_1$ into two solid tori 
$U_1$, $U_2$, 
where $K$ is a core circle of $U_1$. 
If $\partial D_2 \subset \partial U_1$, 
then the above argument works 
to show that 
there exists a non-separating meridian 
disk for $V_1$ giving weak $K$--reducibility 
together with $D_2$. 
If $\partial D_2 \subset \partial U_2$, 
then we take a meridian disk $D_1'$ for $U_1$ 
such that 
$\partial D_1' \subset V_1$, 
and $D_1'$ intersects $K$ transversely in one point. 
We may regard $D_1'$ a 
(non-separating essential) $K$--disk in $V_1$, 
and we have 
$\partial D_1' \cap \partial D_2 = \emptyset$. 
By regarding $D_1'$ as $D_1$, we see that 
we may suppose that $D_1$, $D_2$ are non-separating in 
$V_1$, $V_2$ respectively. 

\medskip
Now we have the following two cases. 

\medskip
\noindent
{\bf Case 1}\qua
$D_1 \cap K = \emptyset$. 

\medskip
Let $T$ be the solid torus obtained from 
$V_1$ by cutting along $D_1$. 
Since $\partial D_2$ is non-separating in $\partial V_2$ 
and $S^3$ does not contain non-separating 2--sphere, 
we see that $\partial D_2$ is an essential simple closed 
curve in $\partial T$. 
Since $S^3$ does not contain non-separating 2--sphere or 
punctured lens spaces, 
$\partial D_2$ is a longitude of $T$, 
and, hence, there is an 
annulus $A$ in $T$ such that 
$\partial A = K \cup \partial D_2$. 
Then $A \cup D_2$ gives a disk bounding $K$, 
and this shows that $K$ is a trivial knot. 

\medskip
\noindent
{\bf Case 2}\qua
$D_1 \cap K \ne \emptyset$. 

\medskip
Let $N = N(D_1;V_1)$, 
$T_1 = c\ell (V_1 - N)$, 
$a_1 = K \cap T_1$, and 
$a_2 = K \cap N$. 
Note that $a_2$ is a core with respect to a natural 
1--handle structure on $N$. 
It is easy to see that $a_1$ is a trivial arc in $T_1$. 
Let $T_2 = V_2 \cup N$. 
We regard $a_2$ as an arc properly embedded in $T_2$. 

\medskip
\noindent
{\bf Claim 2}\qua
$T_2$ is a solid torus and $a_2$ is a trivial arc in $T_2$. 

\medskip
\noindent
{\bf Proof of Claim 2}\qua
Let $T'$ be the solid torus obtained from $V_2$ 
by cutting along $D_2$ and $B' = T' \cup N$. 
By the arguments in Case~1, 
we see that $\partial D_1$ is a longitude of $T'$. 
Hence $B'$ is a 3--ball and $a_2$ is a trivial arc in $B'$. 
Since $V_2$ is obtained from 
$B'$ by identifying two disks in 
$\partial B'$ corresponding to the copies of $D_2$, 
we see that $T_2$ is a solid torus, 
and $a_2$ is a trivial arc in $T_2$. 

Hence we see that $T_1 \cup T_2$ gives a 
genus one 1--bridge position of $K$. 
By the construction of $T_1$, 
we see that $\sigma$ is isotopic to 
an unknotting tunnel associated to 
$T_1 \cup T_2$. 
\end{proof}

\section{Comparing 2--bridge position and an 
unknotting tunnel}

In \cite{R-S1}, Rubinstein-Scharlemann introduced a powerful 
machinery called {\it graphic} 
for studying positions of 
two Heegaard surfaces of a 3--manifold. 
Successively, Dr.~Osamu Saeki and the author 
introduced an orbifold version of their setting, 
and showed that the results similar to Rubinstein-Scharlemann's 
hold in this setting \cite{K-S}. 
In this section, 
we quickly review the arguments and apply it 
to compare decomposing 2--spheres giving 2--bridge positions, 
and genus 2 Heegaard splittings obtained from an unknotting tunnel 
for a 2--bridge knot. 

Let $K$ be a 2--bridge knot, 
that is, 
there exists a genus zero Heegaard splitting 
$B_1 \cup _P B_2$ of 
$S^3$ such that $K \cap B_i$ is a 2--strings trivial arcs 
in $B_i$ $(i=1,2)$. 
Then the unknotting tunnels 
$\tau_1$, $\tau_2$ are contained in 
$B_1$, $B_2$ respectively as in 
Figure~3.1. 

\begin{figure}[ht!]\small
\begin{center}
\leavevmode
\epsfxsize=60mm
\epsfbox{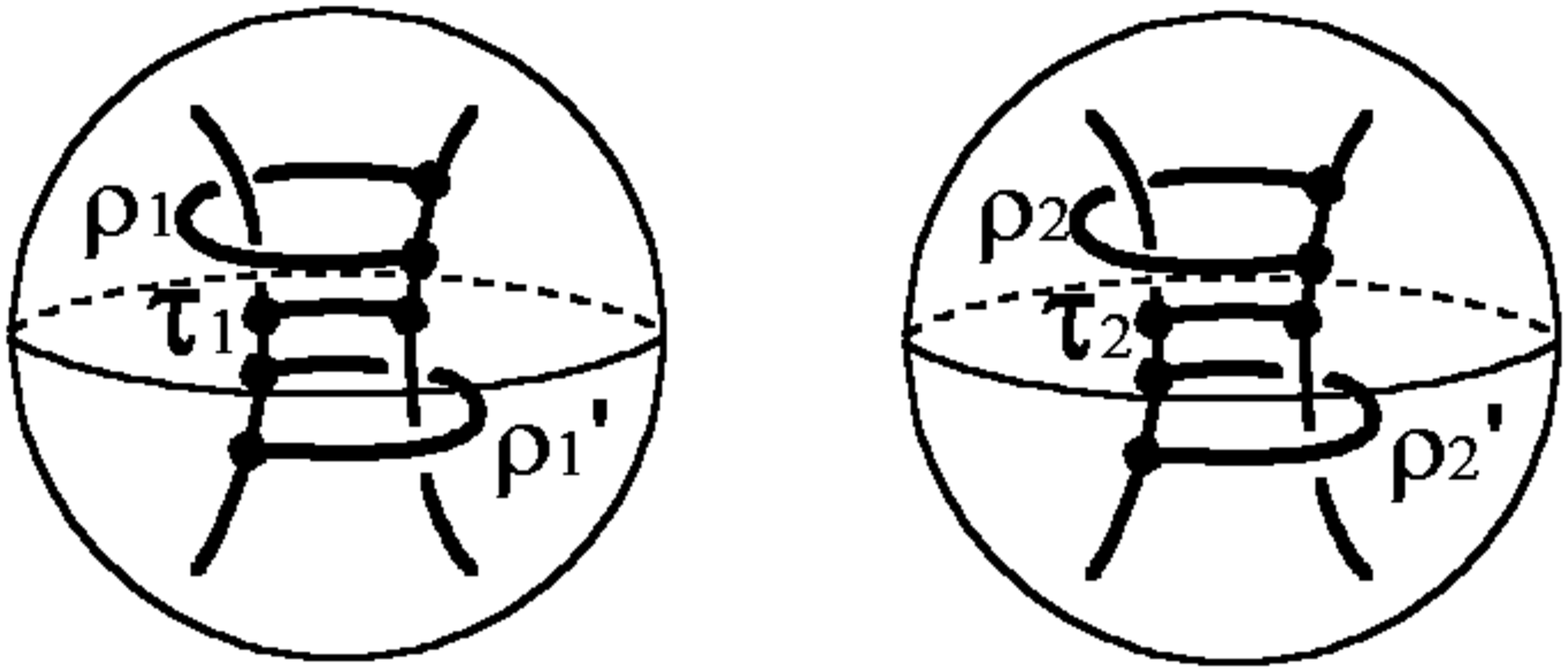}
\end{center}
\begin{center}
Figure 3.1
\end{center}
\end{figure}

There is a diffeomorphism 
$f\co  P \times (0,1) \rightarrow S^3-(\tau_1 \cup \tau_2 )$ 
such that 
$f(P \times \{ 1/2 \})$ is the decomposing 2--sphere $P$, 
and that 
$f((p_1 \cup p_2 \cup p_3 \cup p_4) \times (0,1)) = 
K \cap (S^3-(\tau_1 \cup \tau_2 ))$ 
for some $p_1$, $p_2$, $p_3$, $p_4 \in P$. 

Let $\sigma$ be an unknotting tunnel for $K$. 
Let $\Theta_1 = K \cup \sigma$, 
$V_1 = N(\Theta_1; S^3)$, 
$V_2 = c\ell(S^3-V_1)$, and 
$\Theta_2$ a spine of $V_2$ such that 
each vertex has valency 3. 
Note that 
$V_1 \cup_Q V_2$ is a genus two Heegaard splitting of $S^3$. 
Then there is a diffeomorphism 
$g\co  Q \times (0,1) \rightarrow S^3-(\Theta_1 \cup \Theta_2)$. 

Let 
$P_s = f(P \times \{ s \})$, and 
$Q_t = g(Q \times \{ t \})$. 
Then for a fixed small constant $\varepsilon > 0$, 
we may suppose that 
$P_s \cap Q_t$ looks as one of the following, 
where 
$s \in (0, \varepsilon )$ or $(1-\varepsilon , 1)$, and 
$t \in (0, \varepsilon )$.

\begin{enumerate}

\item 
$P_s \cap Q_t$ consists of two transverse 
simple closed curves $\ell_1$, $\ell_2$ 
which are $K$--essential in $P_s$, 
and inessential in $Q_t$. 

\item 
$P_s \cap Q_t$ consists of a simple closed curve 
$\ell$ 
and a figure $8$ $\delta$ such that; 
$\ell$ is $K$--essential in $P_s$, 
and inessential in $Q_t$, and 
$\delta$ is arising from a saddle tangency. 

\item 
$P_s \cap Q_t$ consists of three transverse 
simple closed curves $\ell_1$, $\ell_2$, and $m$ such that; 
$\ell_1$ and $\ell_2$ bound pairwise disjoint 
$K$--disks in $P_s$ each of which contains a puncture from $K$,
$\ell_1$ and $\ell_2$ are parallel in $Q_t$, and; 
$m$ is $K$--essential in $P_s$ 
and inessential in $Q_t$, 

\item 
$P_s \cap Q_t$ consists of two transverse 
simple closed curves $\ell_1$, $\ell_2$, 
and a figure $8$, $\delta$ 
such that; 
$\ell_1$ and $\ell_2$ bound pairwise disjoint 
$K$--disks in $P_s$ each of which contains a puncture from $K$,
$\ell_1$ and $\ell_2$ are parallel in $Q_t$, and; 
$\delta$ is 
arising from a saddle tangency. 

\item 
$P_s \cap Q_t$ consists of four transverse 
simple closed curves 
$\ell_1$, $\ell_2$, $\ell_3$, and $\ell_4$ 
such that 
$\ell_1$, $\ell_2$, $\ell_3$, $\ell_4$ 
bound mutually disjoint $K$--disks in 
$P_s$ each containing a puncture from $K$, and 
$\ell_1$ and $\ell_2$ 
($\ell_3$ and $\ell_4$ resp.) 
are pairwise parallel in $Q_t$. 

\end{enumerate}

Moreover, 
for a fixed $\varepsilon_1 \in (0, \varepsilon )$, 
if we move $s$ from $0$ to $\varepsilon$, 
then the intersection 
$P_s \cap Q_{\varepsilon_1}$ 
($P_{1-s} \cap Q_{\varepsilon_1}$ resp.) 
is changed as 
$(1) \rightarrow (2) \rightarrow 
(3) \rightarrow (4) \rightarrow (5)$. 

\begin{figure}[ht!]\small
\begin{center}
\leavevmode
\epsfxsize=110mm
\epsfbox{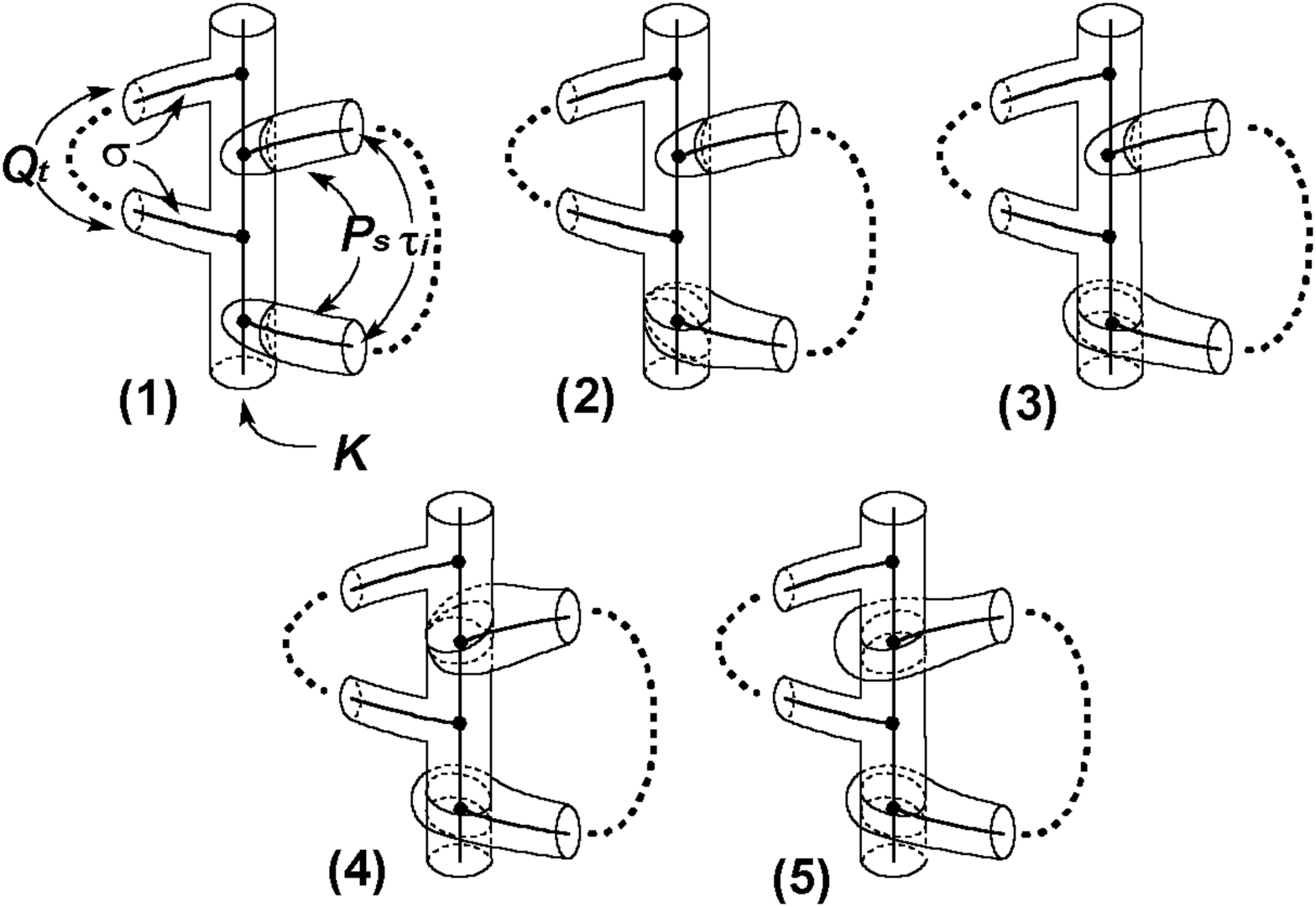}
\end{center}
\begin{center}
Figure 3.2
\end{center}
\end{figure}

Then, by the arguments in Section 4  of \cite{K-S}, 
we see that by an arbitrarily small deformation of 
$f \vert_{(\varepsilon , 1-\varepsilon )}$, and 
$g \vert_{(\varepsilon , 1)}$ 
which does not alter 
$f \vert_{(0, \varepsilon] \cup [1-\varepsilon , 1)}$, and
$g \vert_{(0, \varepsilon ]}$, 
we may suppose that the maps are pairwise generic, 
that is: 

\medskip

There is a stratification of $\text{Int} (I \times I)$ 
which consists of four parts below. 

\begin{description}
\item[Regions] 
Region is a component of the subset of $\text{Int} (I \times I)$ 
consisting of values $(s,t)$ such that 
$P_s$ and $Q_t$ intersect transversely, 
and this is an open set.

\item[Edges]
Edge is a component of the subset consisting 
of values $(s,t)$ such that 
$P_s$ and $Q_t$ intersect transversely except for one 
non-degenerate tangent point. 
The tangent point is either a \lq\lq center\rq\rq \  or a 
\lq\lq saddle\rq\rq. 
Edge is a 1--dimensional subset of $\text{Int}(I \times I)$. 

\item[Crossing vertices]
Crossing vertex is a component of the subset consisting 
of points $(s, t)$ 
such that $P_s$ and $Q_t$ intersect transversely except 
for two non-degenerate tangent points. 
Crossing vertex is an isolated point in $\text{Int}(I \times I)$. 
In a neighborhood of a crossing vertex, four edges are coming in, 
where one can regard the crossing vertex as the intersection of two 
edges. 

\item[Birth-death vertices]
Birth-death vertex is a component of the subset 
consisting of points $(s,t)$ 
such that $P_s$ and $Q_t$ intersect transversely except for 
a single degenerate tangent point. 
In particular, there is a parametrization $(\lambda , \mu )$ 
of $I \times I$ such that 
$P_s = \{ (x,y,z)\vert z=0 \}$, and 
$Q_t =  \{ (x,y,z)\vert z=x^2 +\lambda + \mu y + y^3 \}$. 
Birth-death vertex is an isolated point in $\text{Int}(I \times I)$, 
and
in a neighborhood of a birth-death vertex, 
two edges 
are coming in, 
with one from center tangency, 
the other from saddle tangency. 

\end{description}

Let $\Gamma$ be the union of edges and vertices above. 
By the above, $\Gamma$ is a 1--complex in $\text{Int}(I \times I)$. 
Then we note that as in Section~3 of \cite{R-S1}, 
$\Gamma$ naturally extends to $\partial (I \times I)$. 
Here we note that, 
by the configurations (1) $\sim$ (5) above, 
$\Gamma$ 
looks as in Figure~3.3 near the bottom corners of 
$I \times I$. 
We note that the arguments in Section~6 of \cite{K-S} 
which uses 
labels on the regions hold without changing proofs
in this setting. 
Hence the argument in the proof of 
Proposition~5.9 of \cite{R-S1} which 
uses a simplicial map to a certain 
complex (called $K$ in \cite{R-S1}) 
works in our setting, 
and this shows
(note that $B_1 \cup_P B_2$ is always 
strongly $K$--irreducible 
(Remark of Definition~\ref{2-bridge K-reducible})).

\begin{figure}[ht!]\small
\begin{center}
\leavevmode
\epsfxsize=90mm
\epsfbox{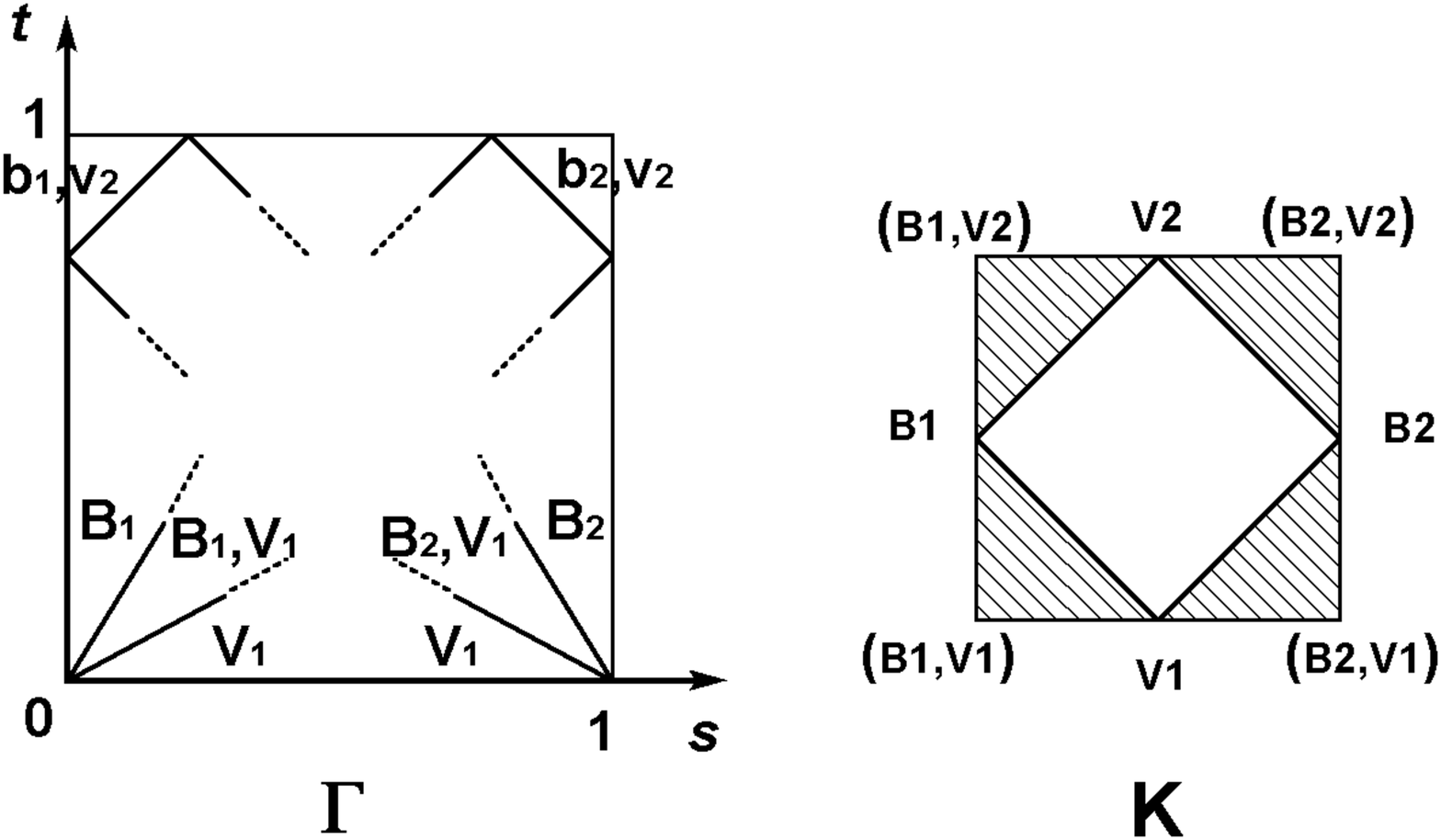}
\end{center}
\begin{center}
Figure 3.3
\end{center}
\end{figure}

\begin{proposition}\label{prop5-6}
Suppose that $V_1 \cup_Q V_2$ 
is strongly $K$--irreducible, 
and 
$K$ is not a trivial knot in $S^3$. 
Then there is an unlabelled region 
in $I \times I - \Gamma$. 
\end{proposition}

And we also have (see Corollary~6.22 of \cite{K-S}): 

\begin{corollary}\label{corollary:essential intersection}
Suppose that 
$V_1 \cup_Q V_2$ 
is strongly $K$--irreducible 
and $K$ is not a trivial knot in $S^3$. 
Then, by applying $K$--isotopy, we may 
suppose that $P$ and $Q$ intersect in non-empty 
collection of simple closed curves which 
are $K$--essential in $P$, 
and essential in $Q$. 
\end{corollary}

\section{Proof of Theorem~\ref{main}}

In this section, we give a proof of  Theorem~\ref{main}. 
For the statements and 
the proofs of Lemmas~B-1, C-1, C-2, C-3, 
D-2, D-3, D-4 which are used in this section, 
see Appendix of this paper. 
Let $K$ be a non-trivial 2--bridge knot and 
$\tau_1$, $\tau_2$, 
$\rho_1$, $\rho_1'$, 
$\rho_2$, $\rho_2'$, 
$\sigma$, 
$B_1 \cup_P B_2$, $V_1 \cup_Q V_2$ 
be as in the previous section. 

\begin{proposition}\label{theorem2}
Suppose that $P \cap Q$ consists of non-empty collection of 
transverse simple closed curves 
which are $K$--essential in $P$ and essential in $Q$. 
Then either 

\begin{enumerate}

\item $\sigma$ is isotopic to either 
$\tau_1$, or $\tau_2$, 

\item $V_1 \cup_Q V_2$ is weakly $K$--reducible, or 

\item there is an essential annulus in $E(K)$. 

\end{enumerate}

\end{proposition}

We note that the closures of $P-Q$ consist of 
two disks 
with each intersecting $K$ in two points, 
and annuli. 
Since the disks are contained in $V_1$, 
$P \cap Q$ consists of even number 
of components. 
The proof of Theorem~\ref{theorem2} is carried out by 
the induction on the number of the components. 
As the first step of the induction, we show: 

\begin{lemma}\label{firststep}
Suppose that $P \cap Q$ consists of two simple 
closed curves 
which are $K$--essential in $P$ and essential in $Q$. 
Then we have the conclusion of Proposition~\ref{theorem2}. 
\end{lemma}

\begin{proof}
Let $D_1$, $A$, $D_2$ be the closures of the components 
of $P-(P \cap Q)$ such that 
$D_1$, $D_2$ are disks, and 
$A$ is an annulus. 

\begin{figure}[ht!]\small
\begin{center}
\leavevmode
\epsfxsize=40mm
\epsfbox{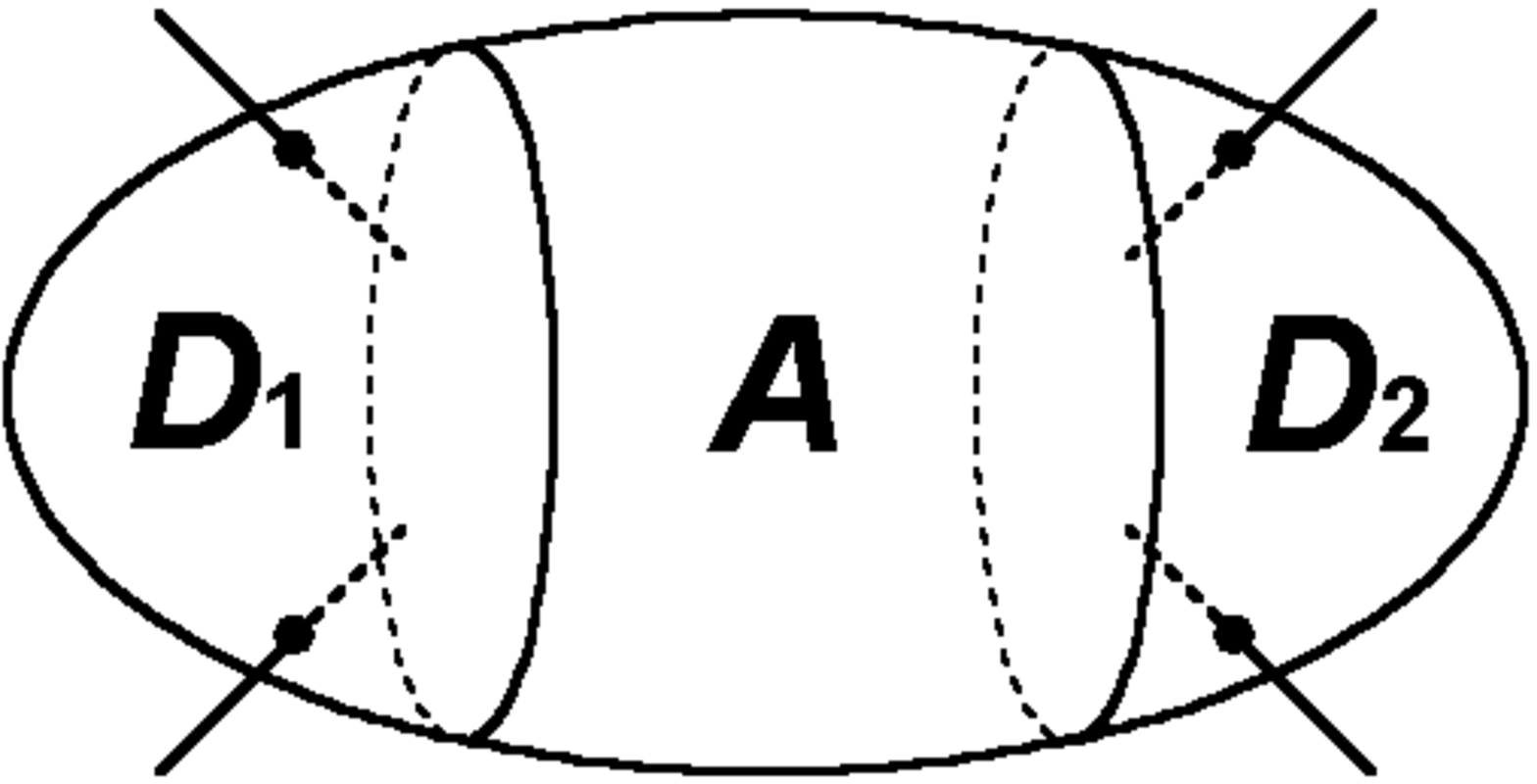}
\end{center}
\begin{center}
Figure 4.1
\end{center}
\end{figure}

We divide the proof into several cases. 

\medskip
\noindent
{\bf Case 1}\qua
Either $D_1$ or $D_2$, say $D_1$, is separating in $V_1$. 

\medskip
We first show: 

\medskip
\noindent
{\bf Claim 1}\qua
The annulus $A$ is boundary parallel in $V_2$. 

\medskip
\noindent
{\bf Proof}\qua
Since $D_1$ is separating in $V_1$, 
the component of $\partial A$ corresponding to 
$\partial D_1$ is 
separating in $\partial V_2$. 
Hence, by Lemma~C-2, we see that $A$ is compressible 
or boundary parallel in $V_2$. 
Suppose that $A$ is compressible in $V_2$. 
Since $S^3$ does not contain non-separating 2--sphere, 
we see that $D_2$ is also separating in $V_1$, 
and, hence, 
$D_1$ and $D_2$ are pairwise parallel in $V_1$. 
Let $A'$ be the annulus in $Q$ such that 
$\partial A' = \partial A$. 
By exchanging suffix, if necessary, 
we may suppose that $A'$ is properly embedded in $B_1$. 
Since each component of $K \cap B_1$ is an unknotted 
arc, we see that $A'$ is an unknotted annulus in $B_1$, 
and this implies that $A$ and $A'$ are parallel in $B_1$, 
and, hence, in $V_2$ 
ie, 
$A$ is boundary parallel. 

This  completes the proof of Claim~1. 

\medskip
By Claim~1, 
we may suppose, by isotopy, that 
$B_1 \subset V_1$, and $\partial B_1 = D_1 \cup A' \cup D_2$, 
where $A'$ is an annulus contained in 
$\partial V_1$ $(= Q)$. 

\medskip
\noindent
{\bf Claim 2}\qua
Both $D_1$ and $D_2$ are $K$--incompressible in $V_1$. 

\medskip
\noindent
{\bf Proof}\qua
Assume, 
without loss of generality, 
that there is a $K$--compressing disk 
$E_1$ for $D_1$. 
Note that since $K \cap D_1$ consists of two points,
$\partial E_1$ and $\partial D_1$ are parallel in $D_1-K$. 
Let $A_1$ be the annulus in $D_1$ bounded by 
$\partial E_1\cup \partial D_1$. 
Let $D_1'$ be the disk in $D_1$ bounded by 
$\partial E_1$. 
Then we have the following two cases. 

\medskip
\noindent
{\bf Case (a)}\qua
$N(\partial E_1; E_1)$ is contained in $B_1$. 

\medskip
We consider the 2--sphere $D_1' \cup E_1$ in $V_1$. 
Let $B_1'$ be the 3--ball in $V_1$ bounded by 
$D_1' \cup E_1$. 
Since $K$ does not contain a local knot in $V_1$, 
we see that $K \cap B_1'$ is an unknotted arc 
properly embedded in $B_1'$. 
Hence there is an ambient isotopy of $S^3$ 
which moves $K \cap B_1'$ to an arc in $D_1$ 
joining $\partial (K \cap B_1')$, 
and which does not move 
$c\ell (K - B_1')$. 
On the other hand, $c\ell(K-B_1')$ is a component 
of the strings of the trivial tangle $(B_2,K \cap B_2)$. 
This shows that $K$ is a trivial knot, 
a contradiction. 

\medskip
\noindent
{\bf Case (b)}\qua
$N(\partial E_1; E_1)$ is contained in $B_2$. 

\medskip
In this case, we first consider the disk 
$A' \cup A_1 \cup E_1$. 
By a slight deformation of $A' \cup A_1 \cup E_1$, 
we obtain a $K$--compressing disk $E_2$ for $D_2$ 
such that 
$N(\partial E_2; E_2)$ 
is contained in $B_1$. 
Then, by the argument as in Case~(a), 
we see that $K$ is a trivial knot, 
a contradiction. 

This completes the proof of Claim~2. 

\medskip
Now we have the following two subcases. 

\medskip
\noindent
{\bf Case 1.1}\qua
$D_1$ and $D_2$ are not $K$--parallel in $V_1$. 

\medskip
In this case, by Lemma~D-4, we see that
$\partial N((K \cup \tau_1); V_1)$ is isotopic to 
$\partial V_1$ in $S^3 - K$. 
This shows that 
$\sigma$ is isotopic to $\tau_1$. 

\begin{figure}[ht!]\small
\begin{center}
\leavevmode
\epsfxsize=50mm
\epsfbox{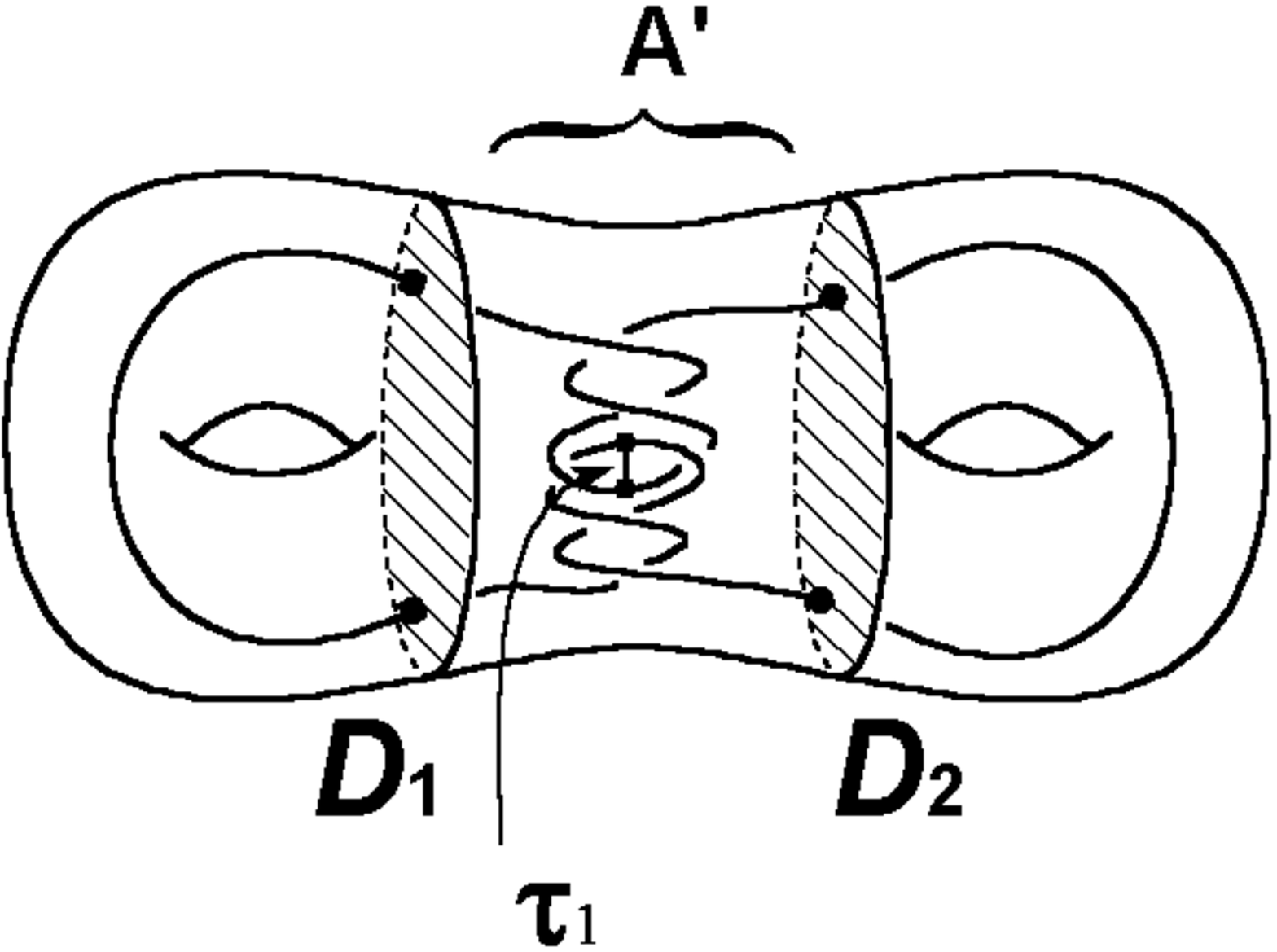}
\end{center}
\begin{center}
Figure 4.2
\end{center}
\end{figure}

\medskip
\noindent
{\bf Case 1.2}\qua
$D_1$ and $D_2$ are $K$--parallel in $V_1$. 

\begin{figure}[ht!]\small
\begin{center}
\leavevmode
\epsfxsize=50mm
\epsfbox{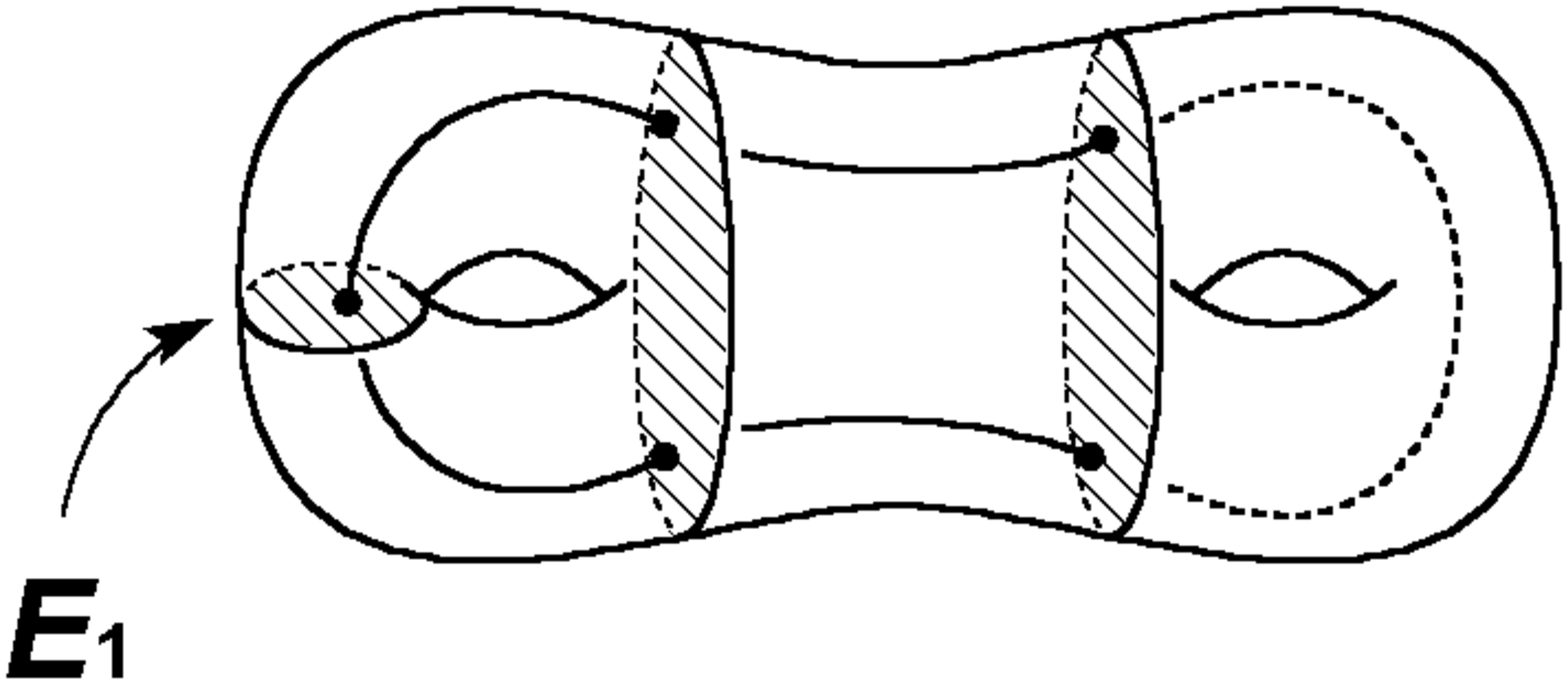}
\end{center}
\begin{center}
Figure 4.3
\end{center}
\end{figure}

Let $Q_1$, $Q_2$ be the closures of the components of 
$Q-A'$ such that 
$\partial Q_i = \partial D_i$ $(i=1,2)$. 
Then $Q_i$ is a torus with one hole 
properly embedded in $B_2$. 
By Lemma~D-2, 
we may suppose, 
by exchanging suffix if necessary, 
that there is a $K$--compressing disk 
$E_1$ for $Q_1$ such that 
$E_1 \subset V_1$, and $E_1 \cap K$ consists of a 
point. 
We consider the genus one surface $Q_2$ 
properly embedded in $B_2$. 
By Lemma~B-1, 
we see that $Q_2$ is $K$--compressible in $B_2$. 
Let $E_2$ be the $K$--compressing disk for $Q_2$. 
Now we have the following subsubcases. 

\medskip
\noindent
{\bf Case 1.2.1}\qua
$N(\partial E_2; E_2)$ is contained in $V_1$. 

\medskip
By the $K$--incompressibility of $D_2$ (Claim~2), 
we see that 
$E_2 \cap K \ne \emptyset$ ie, 
$E_2 \cap K$ consists of a point. 
Then $E_1 \cup E_2$ cuts 
$(V_1, K)$ into a 2--string trivial tangle which is 
$K$--isotopic to $(B_1, K \cap B_1)$. 
Hence $\sigma$ is isotopic to $\tau_1$. 

\begin{figure}[ht!]\small
\begin{center}
\leavevmode
\epsfxsize=55mm
\epsfbox{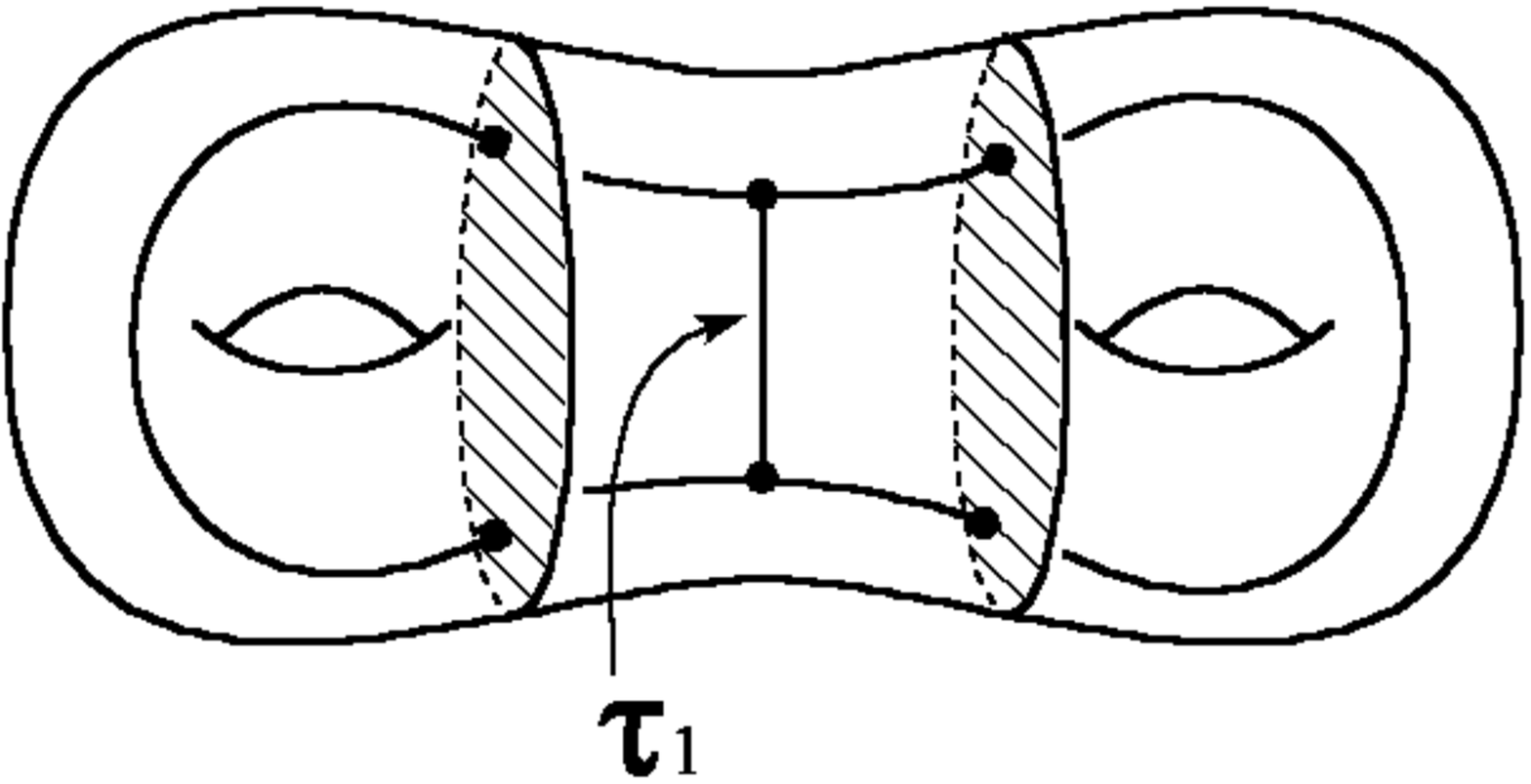}
\end{center}
\begin{center}
Figure 4.4
\end{center}
\end{figure}

\medskip
\noindent
{\bf Case 1.2.2}\qua
$N(\partial E_2; E_2)$ is contained in $V_2$. 

\medskip 
In this case, 
we first show: 

\medskip
\noindent
{\bf Claim 1}\qua
$E_2 \cap Q_1 \ne \emptyset$. 

\medskip
\noindent
{\bf Proof}\qua
Suppose that 
$E_2 \cap Q_1 = \emptyset$. 
Then, by compressing $Q_2$ along $E_2$, 
we obtain a disk $D'$ properly embedded in $B_2$ 
such that $\partial D' = \partial Q_2$, 
and $D'$ separates the components of 
$B_2 \cap K$. 
Let $B_{2,1}$, $B_{2,2}$ be the closures of the 
components of $B_2 - D'$ such that 
$D_1 \subset B_{2,1}$, $D_2 \subset B_{2,2}$. 
Then we can isotope $K \cap B_{2,i}$ rel ~$\partial$ 
in $B_{2,i}$ to an arc in $D_i$ without moving $K \cap B_1$. 
Since $D_1$ and $D_2$ are $K$--parallel in $V_1$, 
this shows that $K$ is a trivial knot, 
a contradiction. 

\medskip
Let $V_{1,2}$ be the closure of the component of 
$V_1 - D_1$ such that 
$\text{Fr}_{B_2} V_{1,2} = Q_1$. 
Note that $V_{1,2}$ is a solid torus in $B_2$ with 
$V_{1,2} \cap P = \partial V_{1,2} \cap P = D_1$. 
By regarding $V_{1,2}$ as a very thin solid torus, 
we may suppose that 
$\text{Int} E_2 \cap V_1$ consists of a 
disk $E_{2,1}$ intersecting 
$K$ in one point. 
Then $E_2 \cap V_2$ is an annulus 
$A_{2,1}$ $(=c\ell (E_2 - E_{2,1}))$. 

\medskip
\noindent
{\bf Claim 2}\qua
$A_{2,1}$ is incompressible in $V_2$. 

\medskip
\noindent
{\bf Proof}\qua
Assume that 
$A_{2,1}$ is compressible in $V_2$.
Then, by compressing $A_{2,1}$, 
we obtain a disk $E_2'$ in $V_2$ 
such that 
$\partial E_2' = \partial E_{2,1}$. 
Since $E_{2,1}$ intersects $K$ in one point, 
$E_{2,1}$ is a non-separating disk in $V_1$.  
Hence, we see that $E_2' \cup E_{2,1}$ 
is a non-separating 2--sphere in $S^3$, 
a contradiction. 

\medskip
Then, by Lemma~C-1, there is an essential disk 
$D_2'$ in $V_2$ such that $D_2' \cap (E_2 \cap V_2) = \emptyset$, 
and, hence, $E_{2,1} \cap D_2' = \emptyset$. 
This shows that $V_1 \cup_Q V_2$ is 
weakly $K$--reducible. 

\begin{figure}[ht!]\small
\begin{center}
\leavevmode
\epsfxsize=55mm
\epsfbox{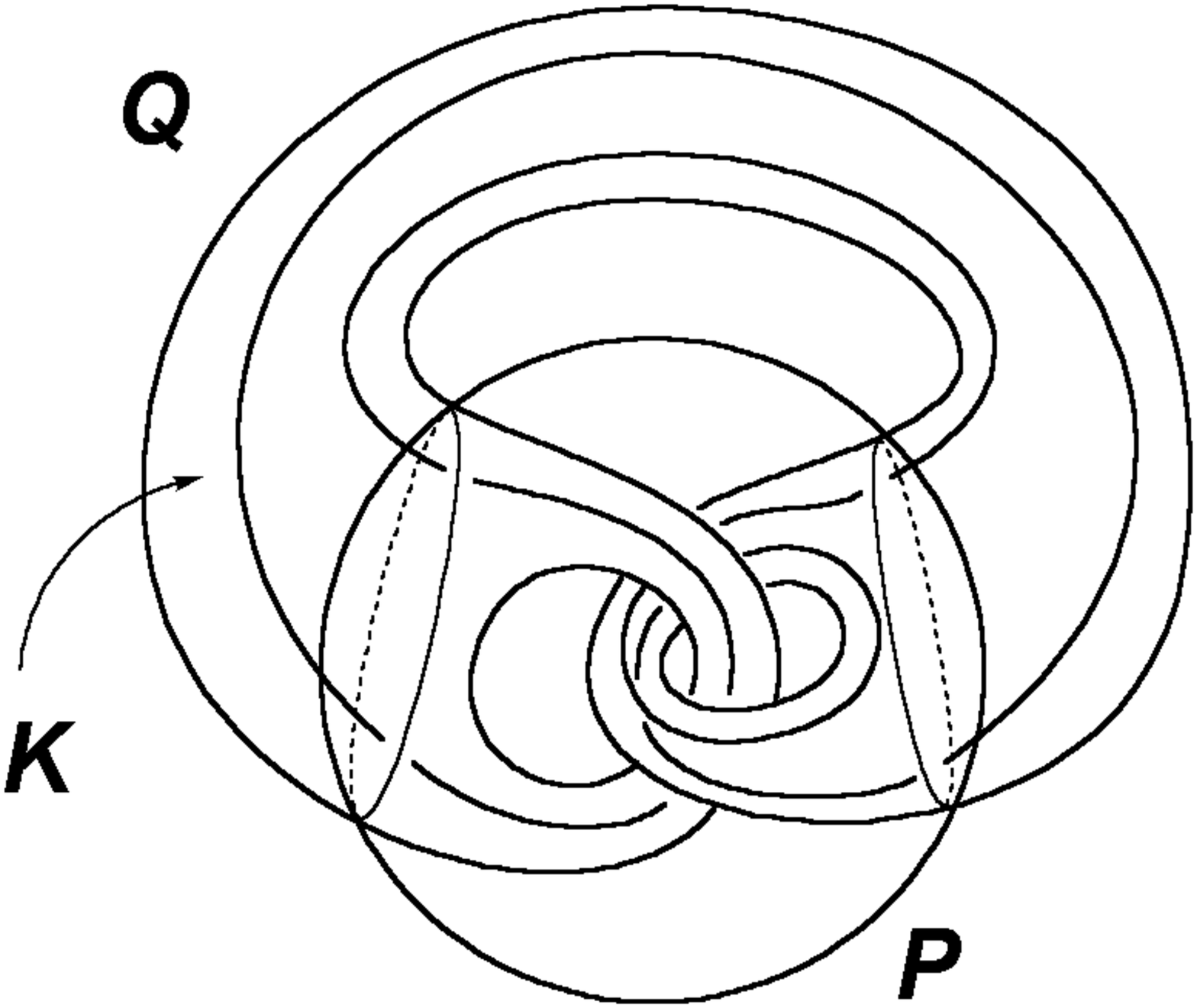}
\end{center}
\begin{center}
Figure 4.5
\end{center}
\end{figure}

\medskip
\noindent
{\bf Case 2}\qua
Both $D_1$ and $D_2$ are non-separating in $V_1$. 

\medskip
In this case, we first show:

\medskip
\noindent
{\bf Claim 1}\qua
$A$ is boundary parallel in $V_2$. 

\medskip
\noindent
{\bf Proof}\qua 
Assume that $A$ is not boundary parallel. 
Since $S^3$ does not contain non-separating 2--sphere, 
we see that $A$ is incompressible in $V_2$. 
Hence, by Lemma~C-1, 
we see that there is an essential disk 
$D$ for $V_2$ such that $D \cap A = \emptyset$, 
and that 
$D$ cuts $V_2$ into two solid tori 
$T_1$, $T_2$, 
where $A \subset T_1$. 
Moreover, since $S^3$ does not contain a punctured 
lens space, we see that each component of 
$\partial A$ represents a generator of the 
fundamental group of the solid torus $T_1$. 
However this contradicts Lemma~C-3.

\medskip
By Claim~1, we may suppose, by isotopy, 
that 
$B_1 \subset V_1$, and 
$\partial B_1 = D_1 \cup A' \cup D_2$, 
where $A'$ is an annulus contained in 
$\partial V_1 (=Q)$. 

\medskip
Then we have the following subcases. 

\medskip
\noindent
{\bf Case 2.1}\qua
Both $D_1$ and $D_2$ are $K$--incompressible in $V_1$. 

\medskip
This case is divided into the following two subsubcases. 

\medskip
\noindent
{\bf Case 2.1.1}\qua
$D_1$ and $D_2$ are not $K$--parallel in $V_1$. 

\medskip
In this case, by Lemma~D-4, 
we see that the given unknotting tunnel $\sigma$ 
is isotopic to $\tau_1$. 

\begin{figure}[ht!]\small
\begin{center}
\leavevmode
\epsfxsize=50mm
\epsfbox{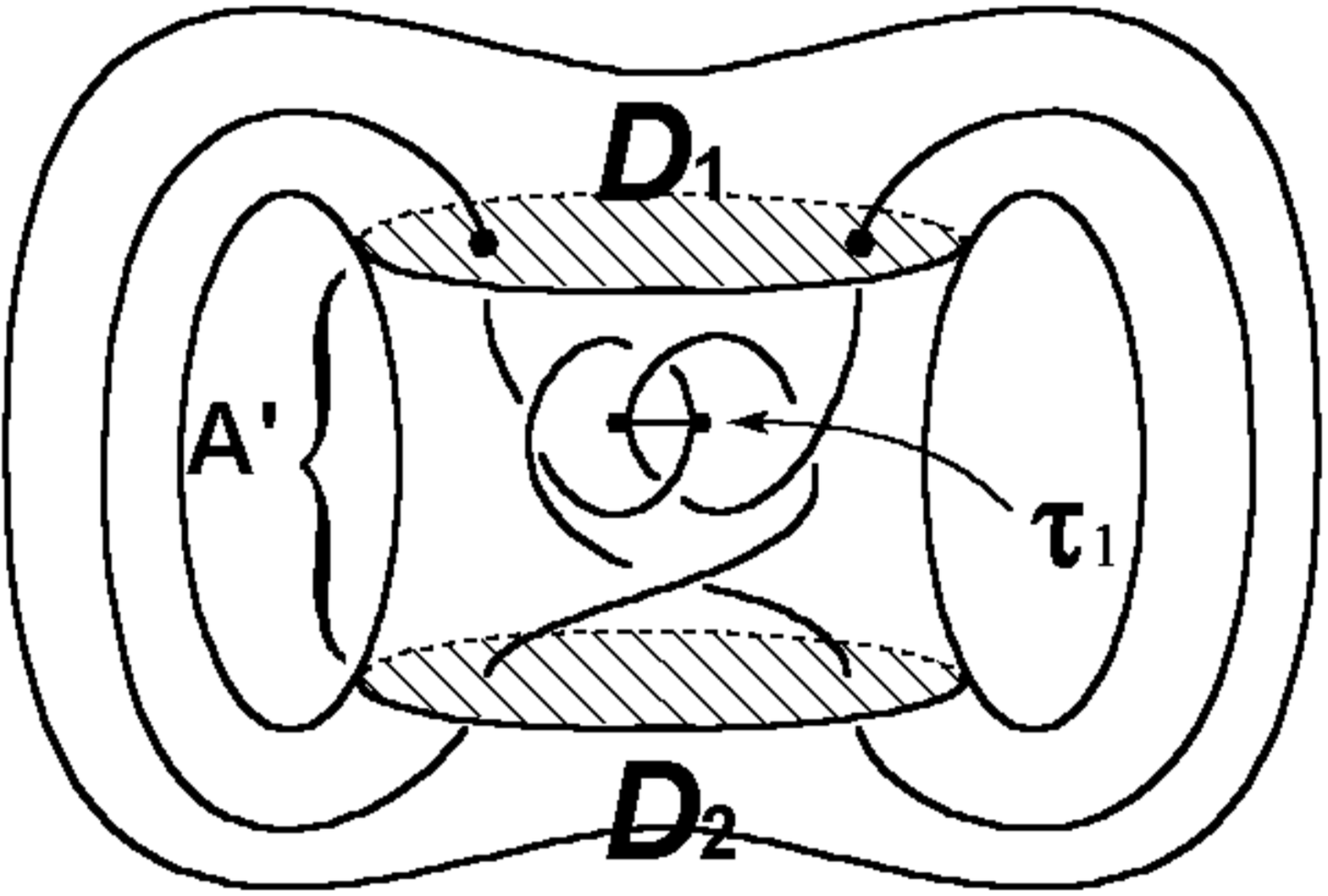}
\end{center}
\begin{center}
Figure 4.6
\end{center}
\end{figure}

\medskip
\noindent
{\bf Case 2.1.2}\qua
$D_1$ and $D_2$ are $K$--parallel in $V_1$. 

\medskip
By Lemma~D-3, there is a $K$--boundary compressing 
disk $\Delta$ for $D_1$ or $D_2$, say $D_1$, 
such that $\Delta \cap D_2 = \emptyset$. 
Let $Q_1$ be the closure of the component of 
$Q-(\partial D_1 \cup \partial D_2)$ 
which is a torus with two holes. 
Let $T_1=Q_1 \cup D_1$. 
Then $\Delta$ is a compressing disk for $T_1$. 
Let $D'$ be the disk obtained by compressing $T_1$ 
along $\Delta$, and 
$D_2'$ a disk obtained by pushing $\text{Int} D'$ 
slightly into $\text{Int} (V_1 \cap B_2)$.  
We may regard $D_2'$ is properly embedded in $B_2$. 
Suppose that $D_2'$ is $K$--compressible in $B_2$. 
Then we can show that $K$ is a trivial knot 
by using the argument as in the proof of 
Claim~1 of Case~1.2.2. 
Hence $D_2'$ is $K$--incompressible in $B_2$. 
Hence, by Lemma~B-1 (3), 
either $D_2'$ and $D_2$ are $K$--parallel 
or $D_2' \cup D_2$ bounds a 2--string trivial tangle in $V_1$, 
which is not a $K$--parallelism between 
$D_2$ and $D_2'$. 

\begin{figure}[ht!]\small
\begin{center}
\leavevmode
\epsfxsize=50mm
\epsfbox{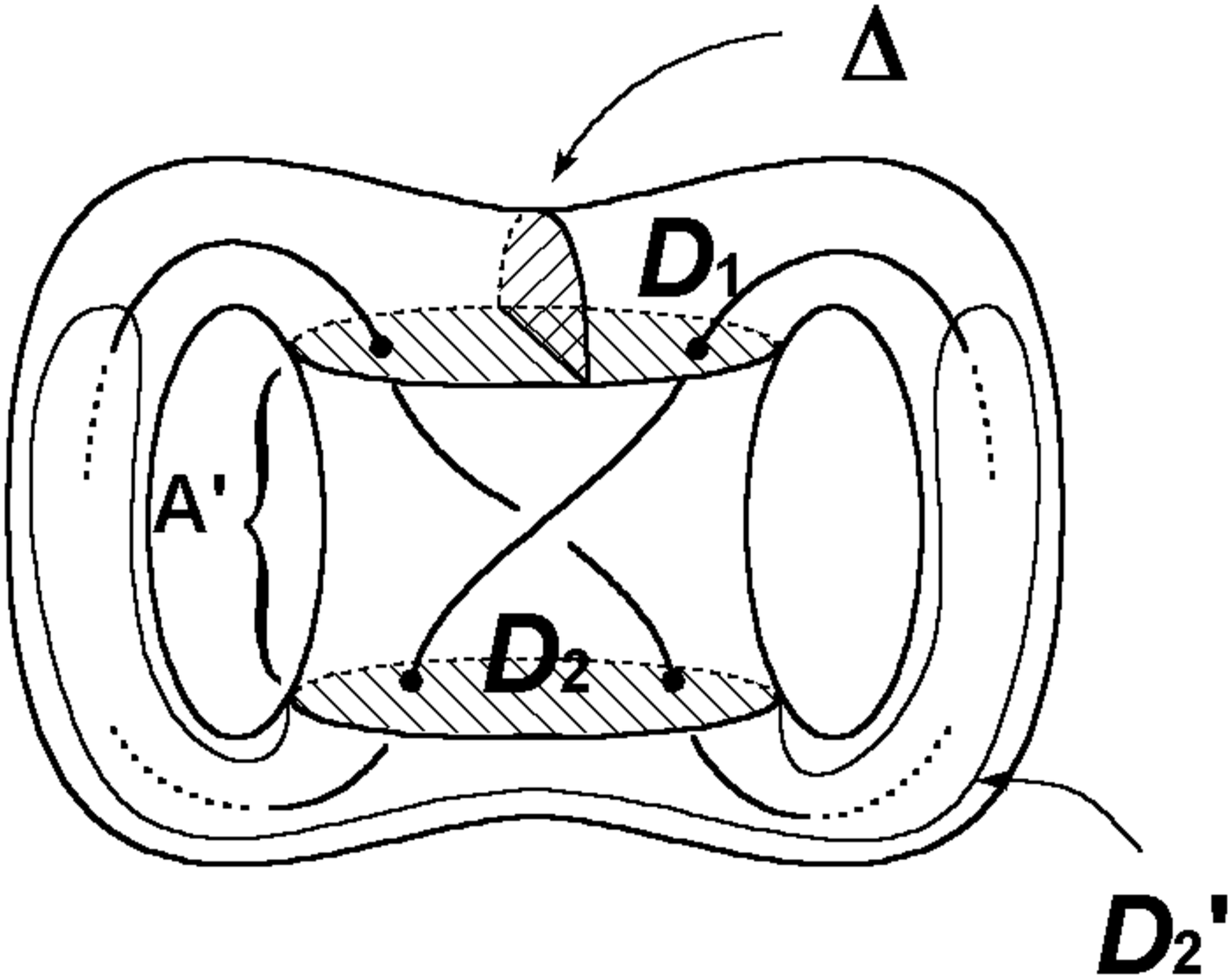}
\end{center}
\begin{center}
Figure 4.7
\end{center}
\end{figure}

In the former case, 
we immediately see that the given unknotting tunnel $\sigma$ is 
isotopic to $\tau_1$. 
In the latter case, 
we have: 

\medskip
\noindent
{\bf Claim 1}\qua
Suppose that $D_2' \cup D_2$ bounds a 
2--string trivial tangle in $V_1$ which is not 
a $K$--parallelism between $D_2$ and $D_2'$. 
Then $\sigma$ is isotopic to $\tau_2$. 

\medskip
\noindent
{\bf Proof}\qua
By Lemma~B-1 (2), 
we see that $D_2'$ and $D_1 \cup A'$ bounds 
a $K$--parallelism in $B_2$. 
Hence, by isotopy, 
we can move $P$ to the position such that 
$B_2 \subset V_1$, and 
$\partial B_2 = D_2 \cup D_2'$. 
Then, by applying the argument of Case~2.1.1 
with regarding $D_2$, $D_2'$ as $D_1$, $D_2$ respectively, 
we see that $\sigma$ is isotopic to $\tau_2$. 

\medskip
\noindent
{\bf Case 2.2}\qua
Either $D_1$ or $D_2$ is $K$--compressible in $V_1$. 

\medskip
Let $E$ be a compressing disk for $D_1$ or $D_2$, 
say $D_1$, in $V_1$. 
Then $\partial E$ and $\partial D_1$ 
are parallel in $D_1-K$, 
and let 
$A^*$ be the annulus in $D_1$ bounded by 
$\partial E \cup \partial D_1$. 
Let $D^*$ be a disk properly embedded in $V_1$ 
which is obtained by moving 
$\text{Int}(A^* \cup E)$ 
slightly so that 
$D^* \cap (D_1 \cup D_2) = 
\partial D^* = \partial D_1$. 

\medskip
\noindent
{\bf Claim 1}\qua
$D^* \subset B_1$. 

\medskip
\noindent
{\bf Proof}\qua
Assume that $D^* \subset B_2$. 
Then we may regard $A' \cup D^*$ is a 
$K$--compressing disk for $D_2$ in $V_1$. 
Then, by using the arguments in Case~(a) of the 
proof of Claim~2 of Case~1, 
we can show that $K$ is a 
trivial knot, a contradiction. 

\medskip
By Claim~1, $\tau_1$ looks as in Figure~4.8.

\begin{figure}[ht!]\small
\begin{center}
\leavevmode
\epsfxsize=45mm
\epsfbox{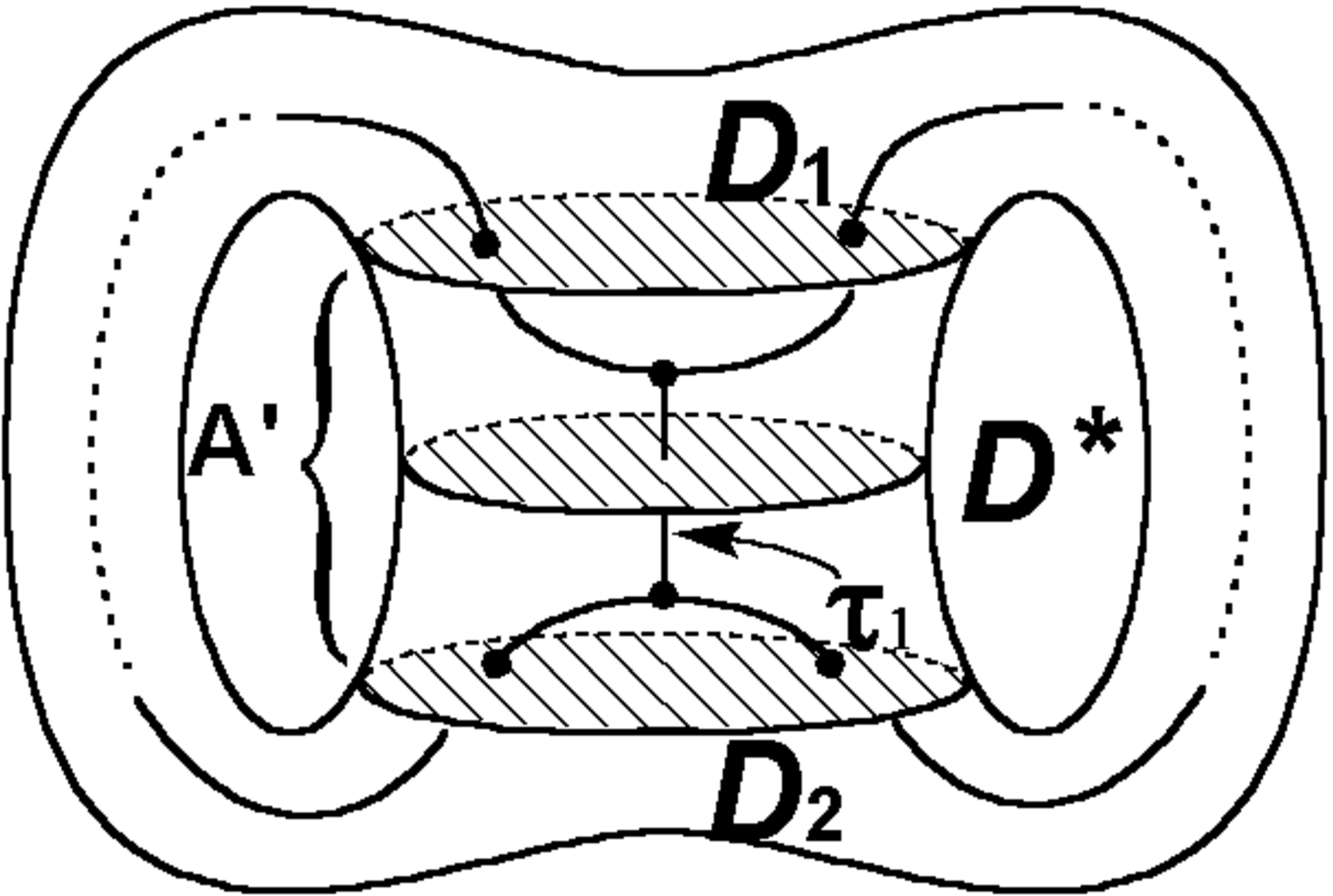}
\end{center}
\begin{center}
Figure 4.8
\end{center}
\end{figure}

%\medskip
\noindent
{\bf Assertion}\qua 
Either 
\lq\lq $K \cup \tau_1$ is a spine of $V_1$\rq\rq  or 
\lq\lq there is an essential annulus in $E(K)$\rq\rq. 

%\medskip
\smallskip
\noindent
{\bf Proof of Assertion}\qua
Let $U_1$ be a sufficiently small regular neighborhood 
of $K \cup \tau_1$, 
and 
$U_2 = c\ell (S^3 - U_1)$. 
Note that $U_2$ is a handlebody, 
because $\tau_1$ is an unknotting tunnel for $K$. 
Let $E_2$ be a non-separating essential disk 
properly embedded $U_2$. 

We may suppose that $D^* \cap U_1$ consists of a 
disk intersecting $\tau_1$ in one point. 

We suppose that $\sharp \{ E_2 \cap D^* \}$ is 
minimal among all non-separating 
essential disks for $U_2$. 

%\medskip
\smallskip
\noindent
{\bf Claim 1}\qua 
No component of $E_2 \cap D^*$ is a simple closed curve, 
an arc joining points in
$\partial U_2$, or 
an arc joining points in $\partial V_1$. 

%\medskip
\smallskip
\noindent
{\bf Proof}\qua
This can be proved by using standard 
innermost disk, outermost arc, and outermost circle 
arguments. 
The idea can be seen in the following figures. 

\begin{figure}[ht!]\small
\begin{center}
\leavevmode
\epsfxsize=100mm
\epsfbox{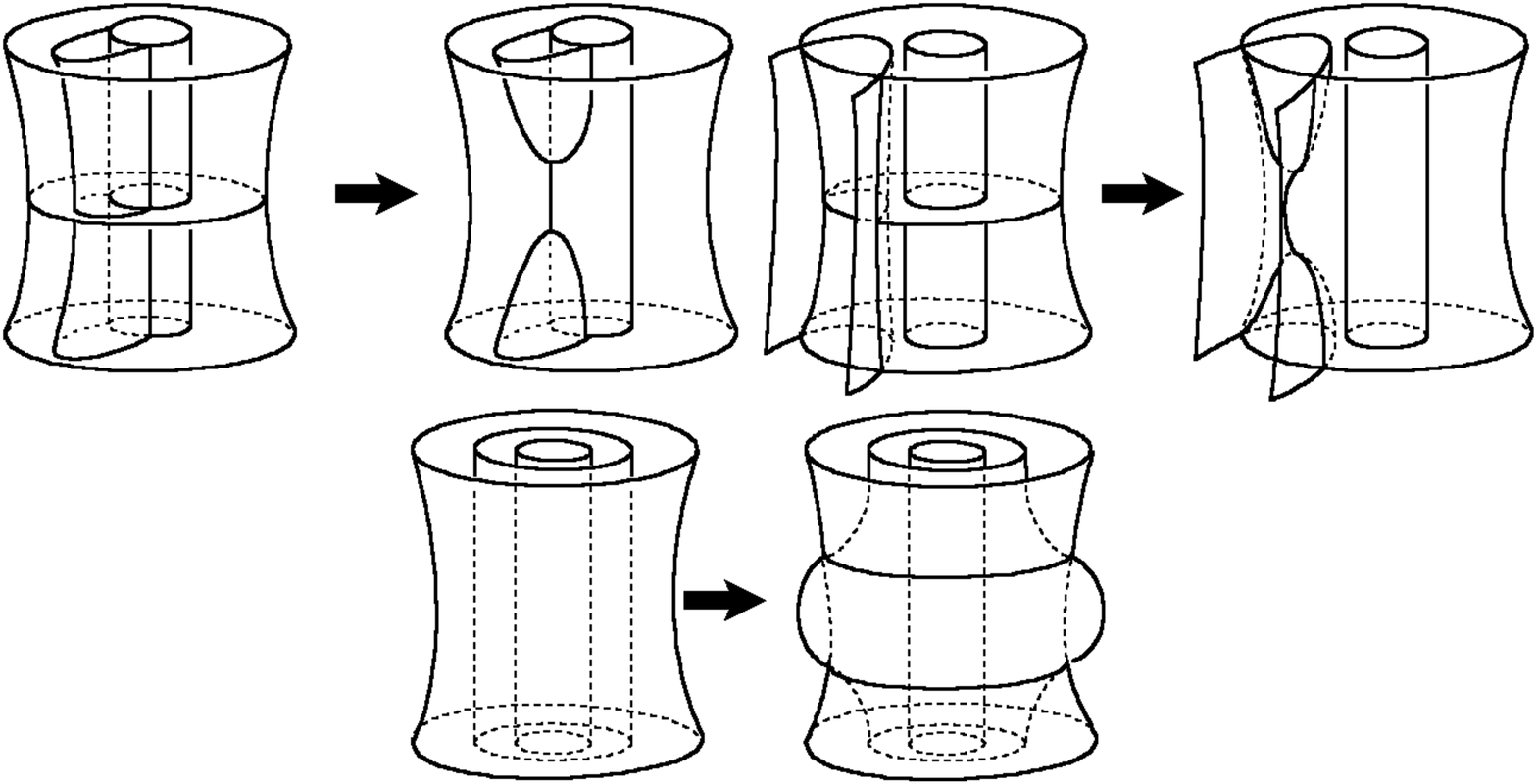}
\end{center}
\begin{center}
Figure 4.9
\end{center}
\end{figure}

\medskip
\noindent
{\bf Claim 2}\qua
$E_2 \cap D^* \ne \emptyset$. 

\medskip
\noindent
{\bf Proof}\qua
Assume that $E_2 \cap D^* = \emptyset$. 
Let $T^*$ be the solid torus obtained by cutting 
$U_1$ along $D^* \cap U_1$. 
Note that $T^*$ is a regular neighborhood of $K$. 
Since $E_2$ is non-separating in $U_2$, and 
$S^3$ does not contain a non-separating 
2--sphere, $\partial E_2$ is an essential simple closed 
curve in $\partial T^*$, and 
$\partial E_2$ is not contractible in $T^*$. 
This shows that $K$ bounds a disk which is 
an extension of $E_2$. 
Hence $K$ is a trivial knot, 
a contradiction. 

\medskip
Hence $E_2 \cap D^*$ consists of a number of arcs 
joining points in 
$\partial U_1$ to points in $\partial V_1$. 
Here, 
by using cut and paste arguments, 
we remove the components of 
$E_2 \cap \partial V_1$ 
which are inessential in $\partial V_1$. 

\medskip
\noindent
{\bf Claim 3}\qua
The components of $E_2 \cap \partial V_1$ 
are not nested in $E_2$. 

\medskip
\noindent
{\bf Proof}\qua
Let $\ell$ be a 
component of $E_2 \cap \partial V_1$ which is 
innermost in $E_2$, and 
$G$ the disk in $E_2$ bounded by $\ell$. 

\medskip
\noindent
{\bf Subclaim 1}\qua
$G$ is contained in $V_2$. 

\medskip
\noindent
{\bf Proof}\qua
Assume that $G$ is contained in $V_1$. 
Since $G \cap (K \cup \tau_1) = \emptyset$, 
this implies that $\tau_1$ is contained in 
a regular neighborhood of $K$, 
contradicting the fact that $\tau_1$ 
is an unknotting tunnel. 

\medskip
\noindent
{\bf Subclaim 2}\qua
$\partial G \cap \partial D^* \ne \emptyset$. 

\medskip
\noindent
{\bf Proof}\qua
Assume that 
$\partial G \cap \partial D^* = \emptyset$. 
Then we can show that there is a non-separating 
disk $G^*$ properly embedded in $V_2$ such that 
$\partial G^* \cap \partial D^* = \emptyset$ 
by using the argument as in the 
Proof of Claim~1 of the proof of 
Proposition~\ref{weakly reducible}. 
Then by using the argument as in the proof of Claim~2 above, 
we can show that $K$ is a trivial knot, 
a contradiction. 

\medskip
Hence there exists a component of 
$E_2 \cap D^*$ connecting 
$\ell$ and $\partial U_1$.
This means that $\ell$ is not surrounded 
by another component of $E_2 \cap \partial V_1$, 
and this gives the conclusion of Claim~3. 

\medskip
\noindent
{\bf Claim 4}\qua
For each component $\ell$ of 
$E_2 \cap \partial V_1$, 
$\ell \cap D^*$ consists of more than one component. 

\medskip
\noindent
{\bf Proof}\qua
Assume that 
$\ell \cap D^*$ consists of a point. 
Let $G$ be the disk in $E_2$ bounded by $\ell$. 
Then $\partial D^*$ and $\partial G$ intersects in 
one point, 
and this shows that $\hat \tau_1$ is a trivial arc 
in $E(K)$, a contradiction. 

\medskip
Let $E^2 = E_2 \cap V_1$. 
We call the boundary component of 
$\partial E^2$ corresponding to 
$\partial E_2$ the {\it outer boundary}. 
Other boundary components of $E^2$ 
(: the components of $E^2 \cap \partial V_1$) are 
called {\it inner boundary components}. 
Let $V_1'$ be the solid torus obtained by cutting 
$V_1$ along $D^*$. 
Let $\ell$ be 
an inner boundary component 
which is \lq\lq outermost\rq\rq \ 
with respect to the intersection 
$E^2 \cap D^*$, that is: 

\medskip
Let $A_\ell$ be the union of the components of 
$E^2 \cap D^*$ intersecting $\ell$. 
Then 
except for at most one component, each component 
of $E^2 - A_\ell$ does not 
intersect other inner boundary components. 

\medskip
Let $G$ be the disk in $E_2$ bounded by $\ell$. 
Let $a_1, \dots , a_n$ be the components of 
$E^2 \cap D^*$, 
which are located on $E^2$ in this order, 
where $a_i \cup a_{i+1}$ $(i=1, \dots , n-1)$ 
cobounds a square $\Delta_i$ in $E^2$. 
Let $\Delta_i' = \Delta_i \cap V_1'$. 

\begin{figure}[ht!]\small
\begin{center}
\leavevmode
\epsfxsize=80mm
\epsfbox{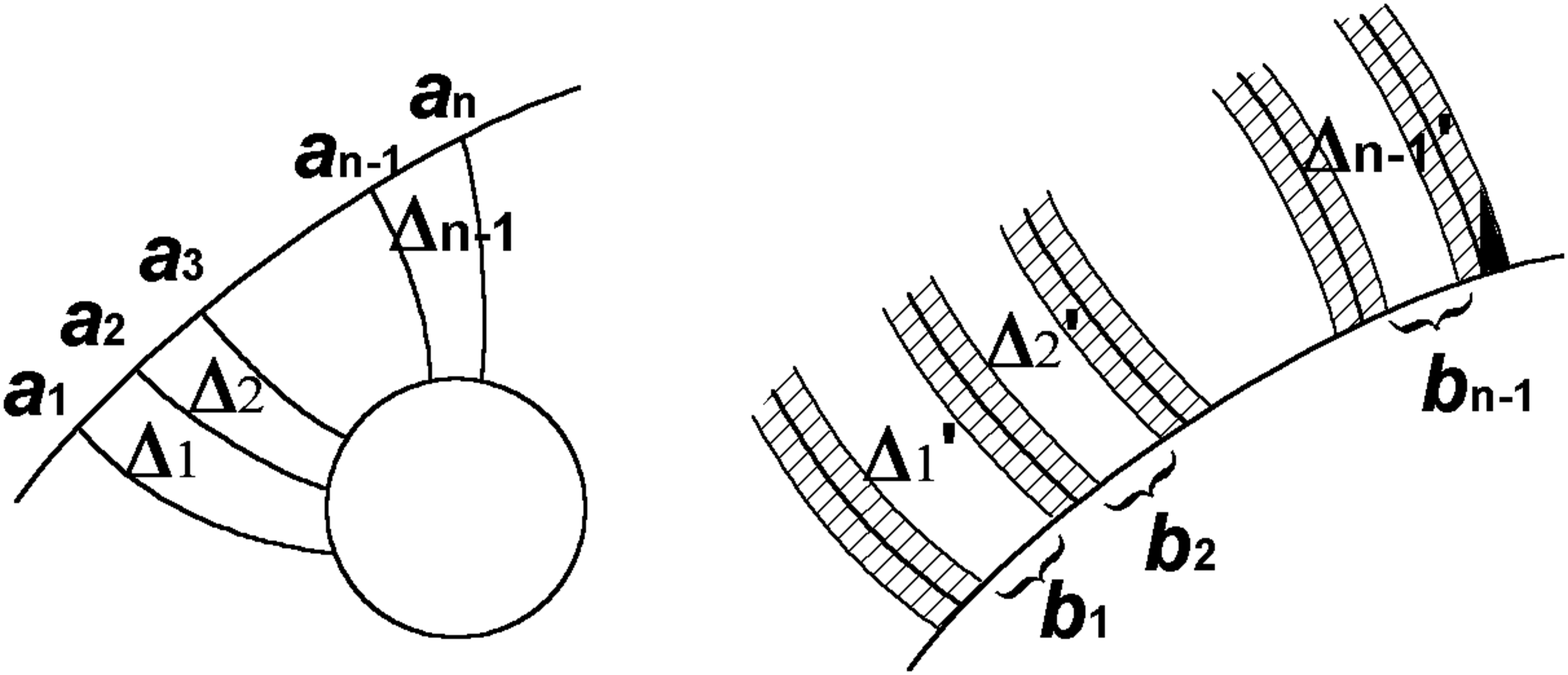}
\end{center}
\begin{center}
Figure 4.10
\end{center}
\end{figure}

Let $R'$ be the image of $\partial V_1$ in $V_1'$. 
Note that $R'$ is a torus with two holes. 
Let $b_i = \Delta_i' \cap R$. 
Then by the minimality condition, we see that 
each $b_i$ is an essential arc properly embedded 
in $R'$. 

\medskip
\noindent
%{\bf Claim 4}\qua
{\bf Claim 5}\qua
If $b_1, \dots , b_{n-1}$ are mutually parallel in $R'$, 
then there is an essential annulus in $E(K)$. 

\medskip
\noindent
{\bf Proof}\qua
Note that $\ell \cap R'$ consists of 
$n$ components, that is, 
$b_1, \dots , b_{n-1}$ above, and 
another component, say $b_0$. 

\medskip
\noindent
{\bf Subclaim 1}\qua
$b_0$ is not parallel to $b_i$ $(i=1, \dots , n-1)$ 
in $R'$. 

\medskip
\noindent
{\bf Proof}\qua
Assume that 
$b_0, b_1, \dots , b_{n-1}$ 
are mutually parallel in $R'$. 
Then we can take a
simple closed curve $m$ in $\partial V_1$ such that 
$m$ intersects 
$\partial D^*$ transversely in one point, and 
$m \cap R'$ is ambient isotopic to $b_i$ in $R'$. 
Let $T^*$ be a regular neighborhood of 
$D^* \cup m$ in $V_1$ such that 
$\partial G \subset T^*$. 
Note that $T^*$ is a solid torus, 
and $\partial G$ wraps around $\partial T^*$ 
longitudally $n$ times. 
This show that the 3--sphere contains a lens 
space with fundamental group a cyclic group of 
order $n$, a contradiction. 

\medskip 
By Subclaim~1, 
we see that we can take simple closed curves 
$m_0$, $m_1$ in $\partial V_1$ 
such that $m_0 \cap m_1 = \emptyset$, 
$m_i$ $(i=0,1)$ intersects $\partial D^*$ 
transversely in one point,  
$m_0 \cap R'$ is ambient isotopic to $b_0$ in $R'$, and 
$m_1 \cap R'$ is ambient isotopic to $b_i$ $(i=1, \dots n-1)$ 
in $R'$.

Let $W^*$ be a regular neighborhood of 
$D^* \cup m_0 \cup m_1$ in $V_1$ 
such that $\partial G \subset \partial W^*$, 
and $A^* = \text{Fr}_{V_1} W^*$. 
Then $W^*$ is a genus two handlebody, and 
$A^*$ is an annulus in $\partial W^*$. 
Note that 
$c\ell(V_1 - W^*)$ is a 
regular neighborhood of $K$. 
Then we denote by $E'(K)$ the closure of 
the exterior of this 
regular neighborhood of $K$. 
Note that 
$A^*$ is embedded in $\partial E'(K)$. 
Then attach $N(G; V_2)$ to 
$W^*$ along $\partial G = \ell$. 
It is directly observed (see Figure~4.11) that 
we obtain a solid torus, say $T_*$, such that $A^*$ 
wraps around $\partial T_*$ longitudally $n$--times. 
Then, let $A^*{}' = c\ell (\partial T_* - A^*)$. 
Note that $A^*{}'$ is an annulus properly embedded 
in $E'(K)$. 

\begin{figure}[ht!]\small
\begin{center}
\leavevmode
\epsfxsize=60mm
\epsfbox{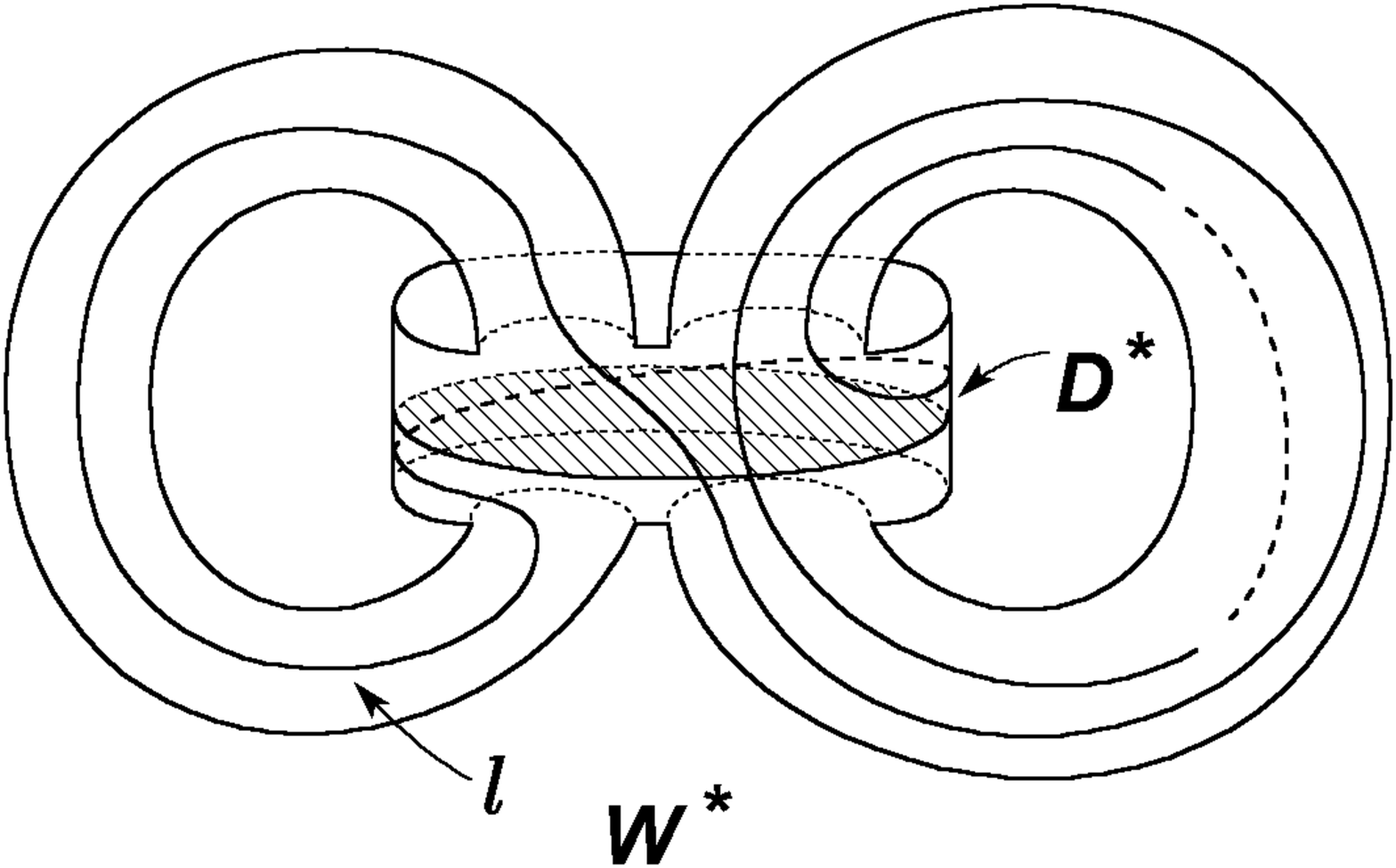}
\end{center}
\begin{center}
Figure 4.11
\end{center}
\end{figure}

Assume that $A^*{}'$ is compressible in $E'(K)$. 
Then the compressing disk is not contained in $T_*$ 
since $A^*{}'$ is incompressible in $T_*$. 
Hence $T_*$ together with a regular neighborhood of 
this compressing disk produces a punctured lens space 
with fundamental group a cyclic group of order $n$ in $S^3$, 
a contradiction. 
Hence $A^*{}'$ is incompressible in $E'(K)$. 
Then assume that $A^*{}'$ is boundary parallel, 
and let $R$ be the corresponding parallelism. 
Since $n \ge 2$, $R$ is not $T_*$. 
Hence $E'(K) = T_* \cup R$, and 
this shows that $E'(K)$ is a solid torus, 
which implies that $K$ is a trivial knot, 
a contradiction. 
Hence $A^*{}'$ is an essential annulus in $E'(K)$, 
and this completes the proof of Claim~5.

\medskip 
Suppose that 
$b_1, \dots , b_{n-1}$ contains at least 
two proper isotopy classes in $R'$. 
We suppose that 
$b_i$, $b_j$ $(i \ne j)$ belong to mutually 
different isotopy classes. 
Let $r_1$, $r_2$ be the components of $\partial R'$. 
Since $\partial G$ and $\partial D^*$ intersects 
transversely, we easily see that we may suppose that 
$b_i \cap r_1 \ne \emptyset$, and 
$b_j \cap r_2 \ne \emptyset$. 

Let 
$T^*$ be the solid torus obtained by cutting $V_1$ along 
$D^*$, and 
$T^2 = c\ell(T^*-N(K;T^*))$ 
(, hence, $T^2$ is homeomorphic to (torus)$\times [0,1]$). 
Here we may regard that $U_1$ is obtained from 
$U_1 \cap T^*$ by adding a 1--handle $h^1$ corresponding to 
$N(D^* \cap U_1; U_1)$, where 
$\tau_1 \cap h^1$ is a core of $h^1$. 
Let $\tau'$, $\tau''$ be the components of the image 
of $\tau_1$ in $T^2$, 
where we may regard that 
$U_1 \cap T^*$ is obtained from 
$N(K, T^*)$ by adding 
$N(\tau' \cup \tau ''; T^2)$. 

\medskip
\noindent
{\bf Claim 6}\qua
$\tau' \cup \tau''$ is \lq\lq vertical\rq\rq \ in
$T^2$ ie, 
$\tau' \cup \tau''$ is ambient isotopic to 
the union of arcs of the form 
$(p_1 \cup p_2) \times [0,1]$, 
where $p_1$, $p_2$ are points in (torus). 

\medskip
\noindent
{\bf Proof}\qua
By extending $\Delta_i'$ ($\Delta_j'$ resp.) to the cores 
of $N(\tau' \cup \tau ''; T^2)$, 
we obtain either an annulus which contains 
$\tau'$  or $\tau''$ 
(if $\partial b_i$ ($\partial b_j$ resp.) is contained 
in $r_1$ or $r_2$), or a rectangle 
two edges of which are $\tau'$ and $\tau''$ 
(if $\partial b_i$ ($\partial b_j$ resp.) joins 
$r_1$ and $r_2$) in $T^2$. 

\begin{figure}[ht!]\small
\begin{center}
\leavevmode
\epsfxsize=90mm
\epsfbox{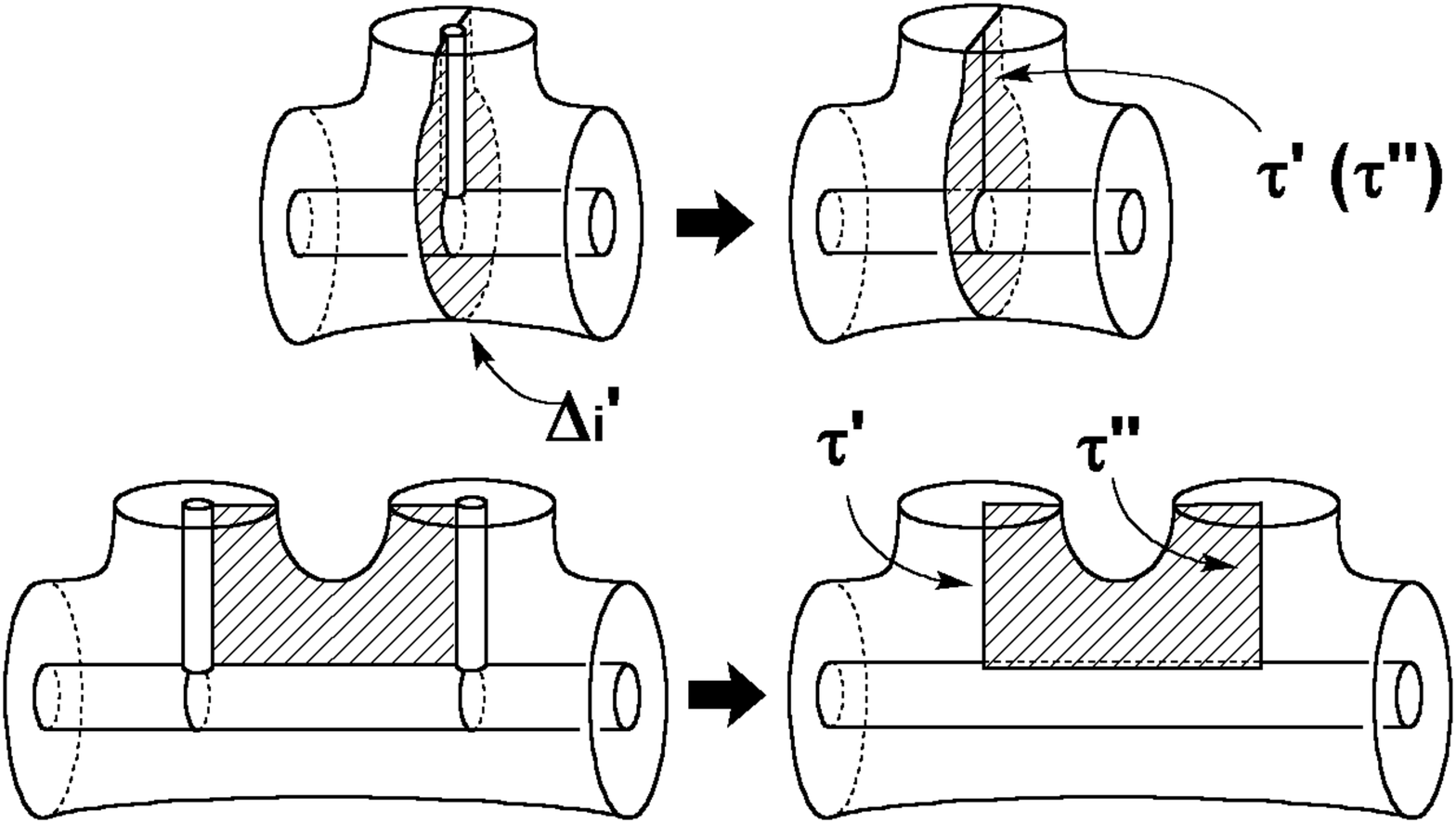}
\end{center}
\begin{center}
Figure 4.12
\end{center}
\end{figure}

Then we have the following three cases. 

\medskip
\noindent
{\bf Case 1}\qua
Both $b_i$ and $b_j$ join $r_1$ and $r_2$. 

\medskip
In this case, we obtain an annulus 
$A_*$ by taking the union of the rectangles 
from $\Delta_i$ and $\Delta_j$. 
Since $b_i$ and $b_j$ are not ambient isotopic in $R'$, 
$A_*$ is incompressible in $T^2$. 
We note that every incompressible annulus in 
$\text{(torus)} \times [0,1]$ with one boundary component 
contained in $\text{(torus)} \times \{ 0 \}$, 
the other in 
$\text{(torus)} \times \{ 1 \}$ is \lq\lq vertical\rq\rq 
(for a proof of this, see, for example, \cite{Ha}). 
Hence $A_*$ is vertical, and this shows 
that $\tau' \cup \tau''$ is vertical. 

\medskip
\noindent
{\bf Case 2}\qua
Either $b_i$ or $b_j$, say $b_i$, join $r_1$ and $r_2$, 
and $\partial b_j$ is contained in $r_1$ or $r_2$. 

\medskip
In this case, we see that $\tau '$ or $\tau''$ is vertical by 
the existence of the annulus from $\Delta_j$. 
Then the existance of the rectangle from $\Delta_i$ 
shows that 
$\tau '$ and $\tau ''$ are parallel, and this implies that 
$\tau' \cup \tau''$ is vertical. 

\medskip
\noindent
{\bf Case 3}\qua
$\partial b_i$ is contained in $r_1$, and 
$\partial b_j$ is contained in  $r_2$. 

\medskip
In this case we see that $\tau' \cup \tau''$ is vertical 
by the existence of the vertical annuli from 
$\Delta_i$ and $\Delta_j$.

\medskip
By Claims~5, and 6, we see that $K \cup \tau_1$ is a spine of 
$V_1$ or there is an essential annulus in $E(K)$, 
and this completes the proof of Assertion. 
\qed

Assertion shows that $\sigma$ is isotopic to $\tau_1$ or 
there is an essential annulus in $E(K)$, 
and this together with the conclusions of 
Cases~1, and 2.1 
shows that we have the conclusions of 
Lemma~\ref{firststep} for all cases. 

This completes the proof of Lemma~\ref{firststep}. 
\end{proof}

\begin{lemma}\label{induction}
Suppose that $P \cap Q$ consists of 
more than two components. 
Then we can deform $Q$ by an ambient 
isotopy in $E(K)$ 
to reduce $\sharp \{ P \cap Q \}$ 
still with non-empty intersection each component of 
which is $K$--essential in $P$, and 
essential in $Q$.  
\end{lemma}

\begin{proof}
Let $2n = \sharp \{ P \cap Q \}$, and 
$D_1, A_1, A_2, \dots , A_{2n-1}, D_2$ 
the closures of the components of 
$P -(P \cap Q)$ such that 
$D_1$, $D_2$ are disks and that they are located 
on $P$ successively in this order. 

\medskip
\noindent
{\bf Claim 1}\qua
Suppose that there is an annulus component 
$A$ of $Q \cap B_i$ ($i=1$ or $2$) such that 
$A$ is $K$--compressible in $B_i$. 
Then the $K$--compressing disk is disjoint from $K$. 

\medskip
\noindent
{\bf Proof}\qua
Let $D$ be the $K$--compressing disk for $A$. 
Assume that $D \cap K \ne \emptyset$ 
ie, 
$D \cap K$ consists of a point. 
Then, by compressing $A$ along $D$, 
we obtain two disks each of which 
intersects $K$ in one point. 
But this is impossible, since each component of 
$\partial A$ separates $\partial B_1$ into two disks 
each intersecting $K$ in two points.

\medskip
\noindent
{\bf Claim 2}\qua
Suppose that there  is an annulus component  
$A_1^Q$ in $Q \cap B_1$, 
and an annulus component $A_2^Q$ in $Q \cap B_2$. 
Then either $A_1^Q$ or $A_2^Q$ is $K$--incompressible 
in $B_1$ or $B_2$. 

\medskip
\noindent
{\bf Proof}\qua
We first suppose that $A_2^Q$ is $K$--compressible in $B_2$. 
Then, by Claim~1, the $K$--compressing disk is disjoint 
from $K$. 
Hence, by compressing $A_2^Q$ along the disk, 
we obtain two disks in $B_2$ which are $K$--essential 
in $B_2$ and disjoint from $K$. 
Let $D_2^*$ be one of the disks. 
Assume, moreover, that 
$A_1^Q$ is also $K$--compressible. 
Then, by using the same argument, 
we obtain a $K$--essential disk $D_1^*$ in $B_1$ 
such that $D_1^* \cap K = \emptyset$. 
Note that $\partial D_1^*$ and $\partial D_2^*$ 
are parallel in $P-K$. 
This implies that $K$ is a two-component 
trivial link, a contradiction. 

\medskip
\noindent
{\bf Claim 3}\qua
If $2n > 6$, then we have the conclusion of 
Lemma~\ref{induction}. 

\medskip
\noindent
{\bf Proof}\qua
Note that there are at most three mutually 
non-parallel, disjoint essential simple closed 
curves on $Q$. 
Hence if $2n > 6$, 
then there are three components, say 
$\ell_1$, $\ell_2$, $\ell_3$, 
of $P \cap Q$ which are mutually parallel on $Q$. 
We may suppose that $\ell_1$, $\ell_2$, $\ell_3$ are
located on $Q$ successively in this order. 
Let $A_1^*$ ($A_2^*$ resp.) be the annulus on $Q$ 
bounded by $\ell_1 \cup \ell_2$ 
($\ell_2 \cup \ell_3$ resp.). 
Without loss of generality, 
we may suppose that $A_1^*$ ($A_2^*$ resp.) is 
properly embedded in $B_1$ ($B_2$ resp.). 
Since $K$ is connected, 
we may suppose, by exchanging suffix if necessary, 
that each component of $\partial A_1^*$ 
separates the boundary points of each 
component of $K \cap B_1$ on $P$. 
Since each component of $K \cap B_1$ is 
an unknotted arc, we see that 
$A_1^*$ is an unknotted annulus. 
Hence there is an annulus $A_1'$ in $P$ 
such that $\partial A_1' = \partial A_1^*$ 
and $A_1'$ and $A_1^*$ are pairwise ($K$--)parallel 
in $B_1$. 
Let $N$ be the parallelism between 
$A_1'$ and $A_1^*$. 

If $\text{Int}(N) \cap Q \ne \emptyset$, 
then we can push the components of 
$\text{Int}(N) \cap Q$ out of $B_1$ 
along the parallelism $N$, still with at least 
two components of intersection $\ell_1 \cup \ell_2$. 
If $\text{Int}(N) \cap Q = \emptyset$, 
then we can push $A_1^*$ out of $B_1$ 
along this parallelism to reduce 
$\sharp \{ P \cap Q \}$ by two. 

\medskip
According to Claim~3 and its proof, 
we suppose that $2n=4$ or $6$, 
and no three components of $P \cap Q$ are mutually 
parallel in $Q$. 
Note that the intersection numbers of any simple closed 
curves on $Q$ with $P \cap Q$ are even, because 
$P$ is a separating surface. 
This shows that $P \cap Q$ consists 
of two (in case when $n=2$) or three (in case when $n=3$) 
parallel classes in $Q$ each of which consists of two 
components. 
Hence, each component of 
$Q \cap B_i$ ($i=1$ or $2$, say 1) is an annulus. 
If a component of $Q \cap B_1$ is $K$--incompressible in 
$B_1$, then, by the argument in the proof of Claim~3, 
we have the conclusion of Lemma~\ref{induction}. 
Hence, 
in the rest of the proof, 
we suppose 
that each component of $Q \cap B_1$ 
is a $K$--compressible annulus in $B_1$. 

Let $N_1$ be the closure of the component of 
$B_1-(Q \cap B_1)$ such that $(K \cap B_1) \subset N_1$. 
Note that $N_1$ is a 3--ball 
such that $\text{Fr}_{B_1}N_1$ consists of 
some components of $Q \cap B_1$. 
Then, by the assumptions, 
we see that $\text{Fr}_{B_1}N_1$ consists of 
either one, two, or three annuli. 

\medskip
\noindent
{\bf Claim 4}\qua
If $\text{Fr}_{B_1}N_1$ consists of an annulus, 
then there is a component of 
$Q \cap B_1$ which is 
$K$--boundary parallel in $B_1$, and, 
hence, we have the conclusion of 
Lemma~\ref{induction}. 

\begin{figure}[ht!]\small
\begin{center}
\leavevmode
\epsfxsize=50mm
\epsfbox{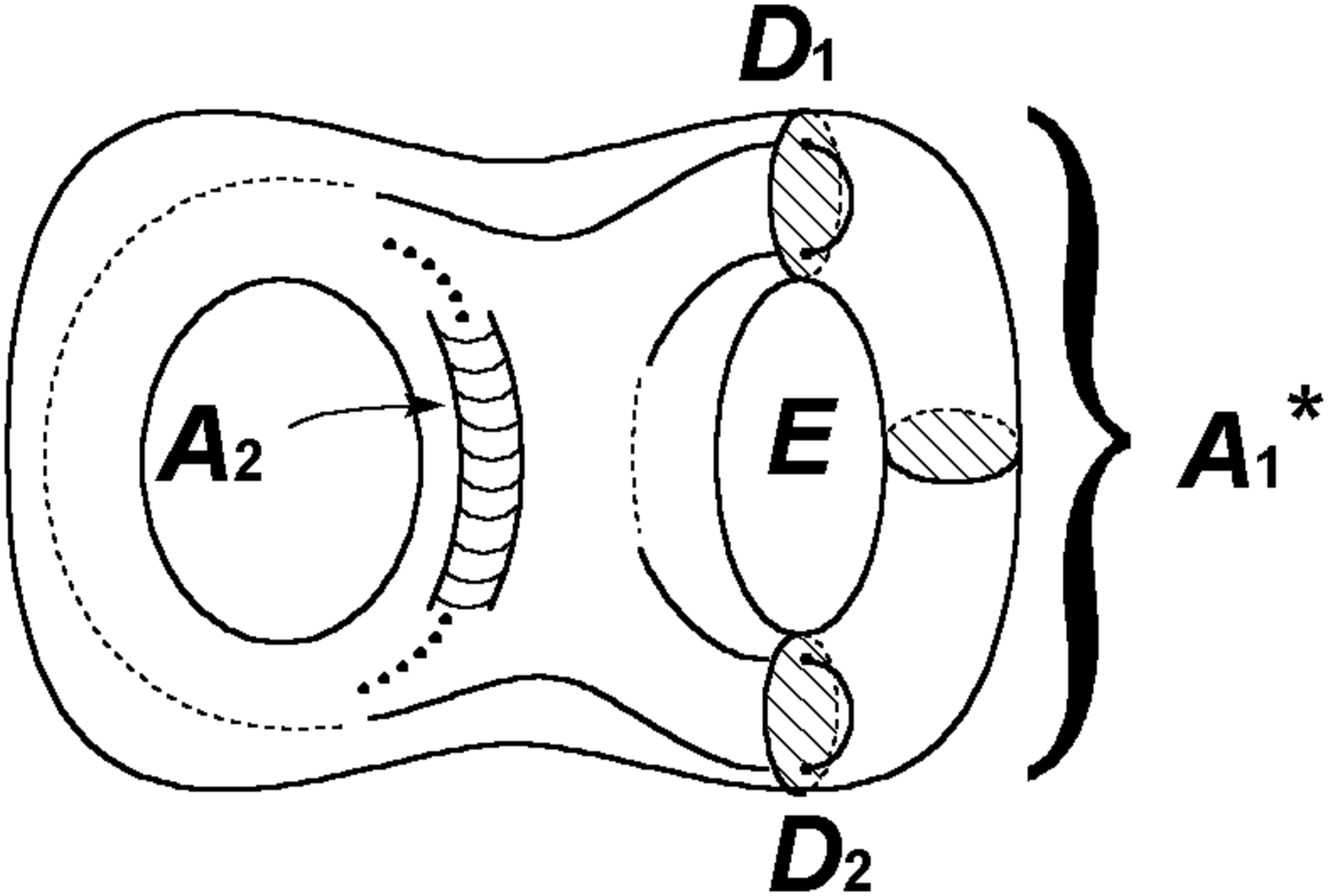}
\end{center}
\begin{center}
Figure 4.13
\end{center}
\end{figure}

\medskip
\noindent
{\bf Proof}\qua
Let $A_1^* = \text{Fr}_{B_1}N_1$. 
Since $A_1^*$ is compressible, 
there is a $K$--compressing disk $E$ for $A_1^*$ 
in $B_1$. 
Note that $E \cap K = \emptyset$ (Claim~1). 
We may regard that 
$E$ is properly embedded in $V_1$ and 
$E$ is parallel to $D_1$ and $D_2$ in $V_1$. 
Since $K$ is connected, 
we see that $E$ is non-separating in $V_1$. 
By cutting $V_1$ along $E$, 
we obtain a solid torus $T_1$ 
such that $K$ is a core circle of $T_1$. 
Recall that 
$D_1, A_1, A_2, \dots , A_{2n-1}, D_2$ 
are the closures of the components of $P-Q$. 
Note that $A_2$ is properly embedded in 
$T_1 - K$. 
Since the 3--sphere does not contain a 
non-separating 2--sphere, 
we see that $A_2$ in incompressible in $T_1$. 
Since every incompressible surface in 
(torus)$\times I$ is either vertical or 
boundary parallel annulus (see \cite{Ha}), 
$A_2$ is boundary parallel in $T_1$. 
Let $N^*$ be the parallelism for $A_2$, 
and $A_2^* = N^* \cap \partial T_1$. 
Since $K$ is connected, and 
$K$ intersects $D_1$ and $D_2$, 
we see that $A_2^*$ is disjoint from 
the images of $D_1$ and $D_2$ in $T_1$. 
Hence we see that 
$A_2^*$ is disjoint from 
the images of $E$ in $\partial T_1$. 
This shows that the parallelism $N^*$ 
survives in $V_1$, and, hence, 
we have the conclusion of Lemma~\ref{induction} 
by the argument as in the proof of Claim~3.

\medskip
\noindent
{\bf Claim 5}\qua
If $\text{Fr}_{B_1}N_1$ consists of two annuli 
$A_1^*$, $A_2^*$, 
then there is a component of 
$Q \cap B_1$ which is 
($K$--)boundary parallel in $B_1$, and, 
hence, we have the conclusion of 
Lemma~\ref{induction}. 

\medskip
\noindent
{\bf Proof}\qua
By exchanging suffix, if necessary, 
we may suppose that the annulus 
$A_i^*$ is incident to $D_i$ $(i=1,2)$. 
If $n=2$, then we have 
$\partial A_1 = \partial A_1^*$. 
If $n=3$, 
then, by reversing the order of 
$A_1, \dots , A_5$, and changing the suffix of 
$A_i^*$ if necessary, we may suppose that 
$\partial A_1 = \partial A_1^*$. 
Then let $N^*$ be the 3--manifold in $B_1$ 
such that 
$\partial N^* = A_1 \cup A_1^*$. 
Note that 
$N^*$ is embedded in $V_2$ and  
$\text{Fr}_{V_2}N^* = A_1$. 

\medskip
\noindent
{\bf Subclaim}\qua
Either $D_1$ or $D_2$, say $D_1$, is 
non-separating in $V_1$. 

\medskip
\noindent
{\bf Proof}\qua
Assume that both $D_1$ and $D_2$ are separating in $V_1$. 
Then $D_1$ and $D_2$ are parallel in 
$V_1$, but this contradicts the fact that 
$N^*$ and $K$ are connected. 

\medskip
Since $D_1$ is a non-separating disk in $V_1$, 
and $S^3$ does not contain a non-separating 
2--sphere, 
we see that $A_1$ is incompressible in $V_2$. 
Then, since $S^3$ does not contain a punctured lens space 
with non-trivial fundamental group, 
we see that 
$A_1$ is boundary parallel in $V_2$ by Lemma~C-3
(see the proof of Claim~1 in Case~2 of 
the proof of Lemma~4.2). 
Hence $N^*$ is a parallelism 
between 
$A_1$ and $A_1^*$, 
and this shows that $A_1^*$ is $K$--boundary 
parallel in $B_1$ along this parallelism 
to give the conclusion of Lemma~\ref{induction}. 

\medskip
\noindent
{\bf Claim 6}\qua
$\text{Fr}_{B_1}N_1$ does not consist of three components. 

\medskip
\noindent
{\bf Proof}\qua
Assume that 
$\text{Fr}_{B_1}N_1$ consists of three annuli 
$A_1^*$, $A_2^*$, and $A_3^*$, 
where 
$\partial A_1^* = \partial A_1$, 
$\partial A_2^* = \partial A_3$, and 
$\partial A_3^* = \partial A_5$. 
Since 
$A_1^*$, $A_2^*$, $A_3^*$ are 
$K$--compressible in $B_1$, 
there are mutually disjoint 
$K$--compressing disks 
$D_1^*$, $D_2^*$, $D_3^*$ for 
$A_1^*$, $A_2^*$, $A_3^*$ respectively. 
We may regard that 
$D_1^*$, $D_2^*$, $D_3^*$ are properly embedded in $V_1$. 
Note that 
$\partial D_1^*$, $\partial D_2^*$, $\partial D_3^*$ 
are not mutually parallel in 
$\partial V_1$. 
Hence we see that 
$D_1^* \cup D_2^* \cup D_3^*$ cuts 
$V_1$ into two components $X_1$, $X_2$ such that 
one component of $K \cap B_1$ is contained in $X_1$, 
and the other component is contained in $X_2$ 
(see Figure~4.14). 
But this contradicts the fact that $K$ is 
connected. 

\begin{figure}[ht!]\small
\begin{center}
\leavevmode
\epsfxsize=60mm
\epsfbox{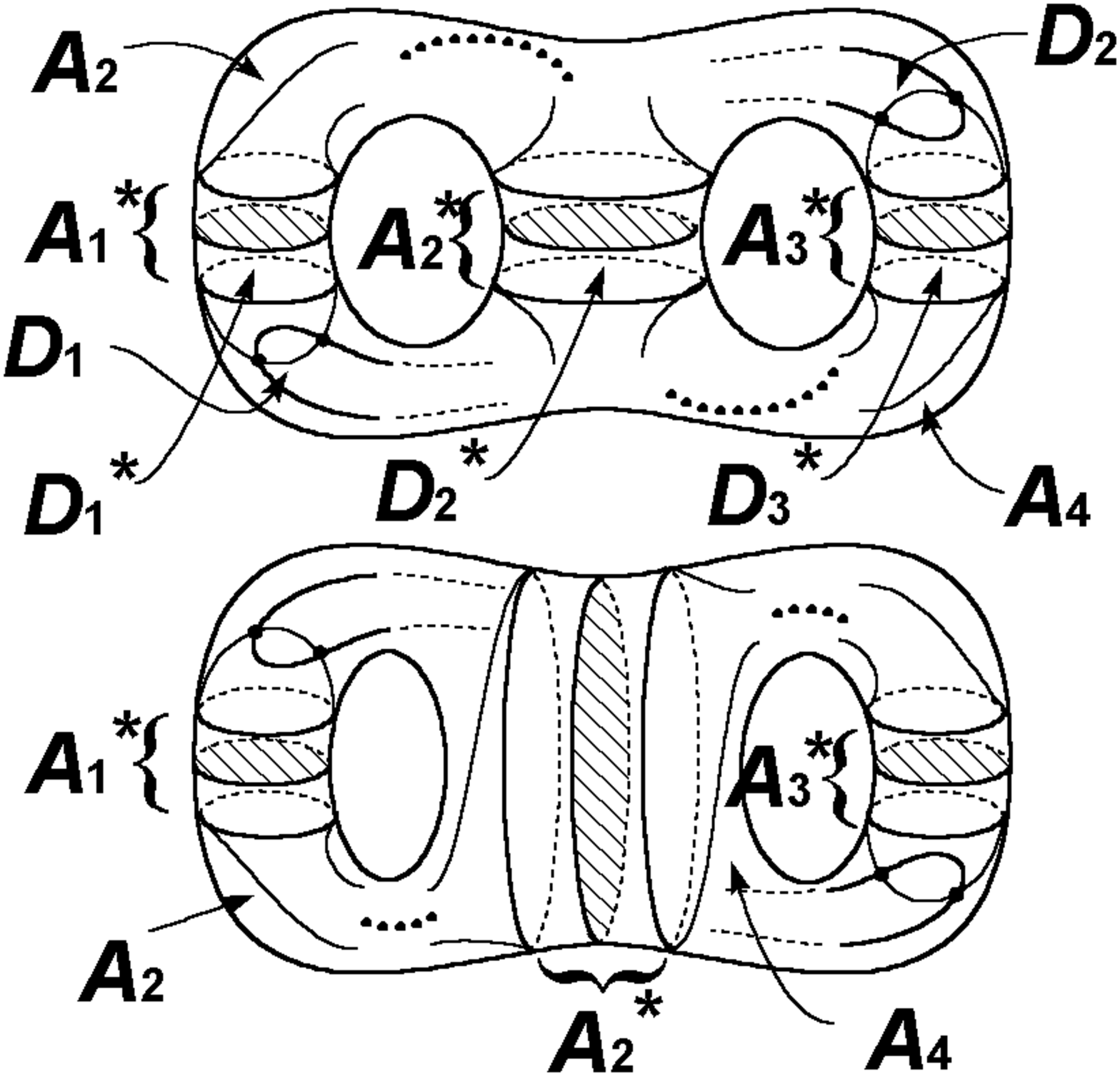}
\end{center}
\begin{center}
Figure 4.14
\end{center}
\end{figure}

\medskip
Claims 3, 4, 5, and 6 complete the proof of 
Lemma~\ref{induction}. 
\end{proof}

\medskip
\noindent
{\bf Proof of Proposition~\ref{theorem2}}\qua
By Lemma~\ref{induction}, 
we may suppose that $P \cap Q$ 
consists of two transverse simple closed curves 
which are $K$--essential in $P$, and 
essential in $Q$. 
Then, by Lemma~\ref{firststep}, 
we have the conclusion of Proposition~\ref{theorem2}. 

\medskip
\medskip
\noindent
{\bf Proof of Theorem~\ref{main}}\qua
Let $\sigma$ be an unknotting tunnel for a non-trivial 
2--bridge knot $K$, and $(V_1, V_2)$ a genus 2 
Heegaard splitting of $S^3$ obtained from 
$K \cup \sigma$ as above. 
If $(V_1, V_2)$ is weakly $K$--reducible, 
then by Propositions~\ref{1to2}, and \ref{weakly reducible}, 
we see that $\sigma$ is isotopic to 
$\tau_1$, $\tau_2$, 
$\rho_1$, $\rho_1'$, 
$\rho_2$, or $\rho_2'$. 
If $(V_1, V_2)$ is strongly $K$--irreducible, 
then by Corollary~\ref{corollary:essential intersection}, 
and Proposition~\ref{theorem2}, 
we see that $\sigma$ is isotopic to 
$\tau_1$ or $\tau_2$, or 
$E(K)$ contains an essential annulus. 
If $E(K)$ contains an essential annulus, 
then $K$ is a $(2, p)$--torus knot. 
Then, by \cite{B-R-S}, it is known that 
every unknotting tunnel for $K$ is isotopic to 
one of $\tau_1$ or $\rho_1$ 
(and that 
$\tau_1$ and $\tau_2$ 
are pairwise isotopic, and 
$\rho_1$, $\rho_1'$, 
$\rho_2$, $\rho_2'$ 
are mutually isotopic). 
Hence we have the conclusion of Theorem~\ref{main}. 

This completes the proof of Theorem~\ref{main}.

%%%%%%%%%%%%%%%%%%%%   End of main body of article
%
%                             References
%

\Addresses\recd

\eject
\small
\setlength{\medskipamount}{4pt plus 2pt minus 2pt}
\setlength{\parskip}{4pt plus 2pt minus 2pt}
\bigskip
{\large\bf Appendix A}

\bigskip
Let $\gamma$ be the union of mutually disjoint 
arcs and simple closed curves properly embedded in a 3--manifold
$N$ such that $N$ admits a 2--fold branched covering space 
$p\co  \tilde{N}\rightarrow N$ along $\gamma$. 

Let $F$ be a surface properly embedded in $N$,
which is in general position with respect to $\gamma$. 
Then, by using ${\Bbb Z}_2$--equivariant loop theorem 
\cite{K-T}, we see that: 

\medskip
\noindent{\bf Lemma A-1}\qua 
$F$ is $\gamma$--incompressible if and only if 
$\tilde{F}$ $(=p^{-1}(F))$ is incompressible.

%Let $a$ be an arc properly embedded in $F$ 
%with $a \cap \gamma = \emptyset$.
%
%\medskip
%\noindent{\bf Definition}\qua 
%We say that $a$ is {\it $\gamma$--inessential} 
%if there is a subarc $b$ of $\partial F$ 
%such that $\partial b = \partial a$, 
%and $a \cup b$ bounds a disk $D$ in $F$ such that 
%$D \cap \gamma = \emptyset$, 
%and $a$ is {\it $\gamma$--essential} if it is not 
%$\gamma$--inessential. 
%
%
%We say that $a$ is {\it weakly $\gamma$--inessential} 
%if there is a subarc $b$ of $\partial F$ 
%such that $\partial b = \partial a$, 
%and $a \cup b$ bounds a $\gamma$--disk $D$ 
%in $F$. 
%
%\medskip
%\noindent{\bf Definition}\qua 
%We say that $F$ is {\it $\gamma$--boundary compressible} 
%if there is a disk $\Delta$ in $N$ such that 
%$\Delta \cap F = \partial \Delta \cap F = \alpha$ 
%is an $\gamma$--essential arc in $F$, and 
%$\Delta \cap \partial N = \partial \Delta \cap \partial N = 
%c\ell (\partial \Delta - \alpha)$. 

\medskip
Moreover, by using ${\Bbb Z}_2$--equivariant cut and paste 
argument as in \cite[Proof of 10.3]{He}, 
we see that: 

\medskip
\noindent{\bf Lemma A-2}\qua 
$\gamma$--incompressible surface $F$ is $\gamma$--boundary compressible
if and only if $\tilde{F}$ is 
boundary compressible. 

\medskip
By using ${\Bbb Z}_2$--Smith conjecture (\cite{Wa}, \cite{Smith}) 
together with the ${\Bbb Z}_2$--equivariant cut and paste 
argument and the irreducibility of $H$, we have: 

\medskip
\noindent{\bf Lemma A-3}\qua 
$\gamma$--incompressible surface $F$ is $\gamma$--boundary parallel 
if and only if $\tilde{F}$ is boundary parallel.
In particular, if $\tilde{N}$ is irreducible, and 
$F$ is a disk intersecting $\gamma$ 
in one point, and $\partial F$ bounds a disk $D$ in 
$\partial H$ such that $D$ intersects $\gamma$ in one 
point, then $F$ is $\gamma$--boundary parallel 
(in fact, $F$ and $D$ are $\gamma$--parallel). 

\bigskip\bigskip 
\bigskip 
{\large\bf Appendix B} 
%($\beta$--essential disks in 2--string trivial tangle)
\bigskip

For the proof of the following two lemmas, 
we refer Appendix~B, and Appendix~C of 
\cite{K-S}.

\medskip
Let $(B, \beta )$ be a 2--string trivial tangle.

\medskip
\noindent{\bf Lemma B-1}\qua 
Let $F$ be a $\beta$--incompressible surface in $B$. 
Then either: 

\begin{enumerate}

\item 
$F$ is a disk disjoint from $\beta$, and 
$F$ separates the components of $\beta$.  
Particularly, in this case, $F$ is $\beta$--essential, 

\item 
$F$ is a $\beta$--boundary parallel disk 
intersecting $\beta$ in at most one point, 

\item 
$F$ is a $\beta$--boundary parallel 
disk intersecting $\beta$ in two points and 
$F$ separates $(B, \beta )$ into the parallelism 
and a rational tangle, or 

\item 
$F$ is a $\beta$--boundary parallel annulus such that 
$F \cap \beta = \emptyset$. 

\end{enumerate}

%\medskip
%\noindent{\it Proof}\qua
%Let $\tilde{F}$ be the lift of $F$ in $\tilde{B}$. 
%By Lemma~A-1, we see that 
%$\tilde{F}$ is one of (1), (2), or (3) of Lemma~B-0. 
%It is easy to see that (1) ((2) resp.) of Lemma~B-0 corresponds 
%to the conclusion (0) ((1) resp.). 
%Suppose that $\tilde{F}$ is an incompressible annulus 
%((3) of Lemma~B-0). 
%Then it is easy to see that we have conclusion (2) if 
%$\tilde{F} \cap \tilde{\beta} \ne \emptyset$, and 
%we have conclusion (3) if 
%$\tilde{F} \cap \tilde{\beta} = \emptyset$. 
%\qed

\medskip
Let $\alpha$ be a 1--string trivial arc in a solid torus $T$. 

\medskip
\noindent{\bf Lemma B-2}\qua 
Let $D$ be an $\alpha$--essential disk in $T$ 
such that $D \cap \alpha$ consists of two points. 
Then there exists an $\alpha$--compressing disk 
$D'$ for $\partial T$ such that 
$D' \cap D = \emptyset$ and 
$D' \cap \alpha$ consists of one point. 
Moreover, by cutting 
$(T, \alpha )$ along $D'$, we obtain a 2--string 
trivial tangle $(B, \beta )$ such that 
$D$ is a $\beta$--incompressible disk in $(B, \beta )$ 
(hence, $D$ is $\beta$--boundary parallel).

\bigskip\bigskip\bigskip 
{\large\bf Appendix C} 

\bigskip
Let $H$ be a genus 2 handlebody, and 
$A$ an essential annulus properly embedded in $H$. 

\medskip
\noindent{\bf Lemma C-1}\qua 
There exists an essential disk 
$D$ in $H$ such that $A \cap D = \emptyset$. 
Moreover the disk $D$ can be taken as a separating disk, or 
a non-separating disk according as 
$A$ is separating or non-separating. 

\medskip
\noindent{\bf Proof}\qua
There exists boundary compressing disk 
$\Delta$ for $A$. 
Apply a boundary compression on $A$ along 
$\Delta$ to obtain a disk $D'$. 
By moving $D'$ by a tiny isotopy, 
we obtain a desired disk $D$. 
For a detail, see, for example, \cite{K2}. 
\qed

\medskip
\noindent{\bf Lemma C-2}\qua 
Each component of $\partial A$ is non-separating in $\partial H$. 
And $A$ is separating in $H$ if and only if the 
components of $\partial A$ are pairwise parallel in $\partial H$. 

\medskip
\noindent{\bf Proof}\qua
Let $D$ be as in Lemma~C-1. 
By the proof of Lemma~C-1, 
we see that $A$ is isotopic to an annulus obtained from 
$D$ by adding a band. 
By isotopy, we may suppose that $A \cap D = \emptyset$. 
Let $T$ be the closure of the component of 
$H-N(D;H)$ such that $A \subset T$. 
Then $T$ is a solid torus, and 
$A$ is incompressible in $T$. 
Hence each component of $\partial A$ is non-separating in 
$\partial T$. 
This implies that each component of 
$\partial A$ is non-separating in $\partial H$. 
Let $D_1$, $D_2$ be the copies of $D$ in $\partial T$. 
Note that $\partial A$ separates $\partial T$ into 
two annuli, say $A_1$, $A_2$. 
If $D$ is separating in $H$, then 
$D_1 \cup D_2$ is contained in one of 
$A_1$ or $A_2$, say $A_1$. 
Then the components of $\partial A$ are 
mutually parallel in $\partial H$ through the 
annulus $A_2$. 
If $D$ is non-separating in $H$, 
then, by exchanging the suffix if necessary, 
we may suppose that $D_1$ is contained in $A_1$, 
and $D_2$ is contained in $A_2$. 
This shows that the components of 
$\partial A$ are not parallel in $\partial H$. 
\qed

\medskip
\noindent{\bf Lemma C-3}\qua 
Let $D$ be as in Lemma~C-1. 
Suppose that $A$ is separating in $H$. 
Then each component of $\partial A$ 
does not represent a generator of the 
fundamental group of the solid torus obtained from 
$H$ by cutting along $D$, which contains $A$. 
%Suppose that $A$ is separating in $H$, and 
%a component of $\partial A$ represents a generator 
%of the fundamental group of the solid torus 
%obtained from $H$ by cutting along $D$. 
%Then $A$ is boundary parallel. 

\medskip
\noindent{\bf Proof}\qua
Let $T$ be the solid torus obtained from $H$ by 
cutting along $D$ such that $A \subset T$. 
Then $(T, A)$ is homeomorphic to 
$(A \times I, A \times \{ 1/2 \})$ 
as pairs. 
This shows that the closure of 
a component of $T - A$ gives a parallelism 
between $A$ and a subsurface of $\partial H$. 
\qed

\eject

{\large\bf Appendix D} 

\bigskip
Let $K$ be a knot in a genus two handlebody $H$ with an essential disk
$E$ such that $E$ cuts $H$ into a solid torus, where $K$ is a core
circle of $T$.  Note that there exists a two-fold branched cover $p\co
\tilde{H}\rightarrow H$ of $H$ along $K$, where $\tilde{H}$ is a genus
three handlebody.

\begin{figure}[ht!]\small
\begin{center}
\leavevmode
\epsfxsize=40mm
\epsfbox{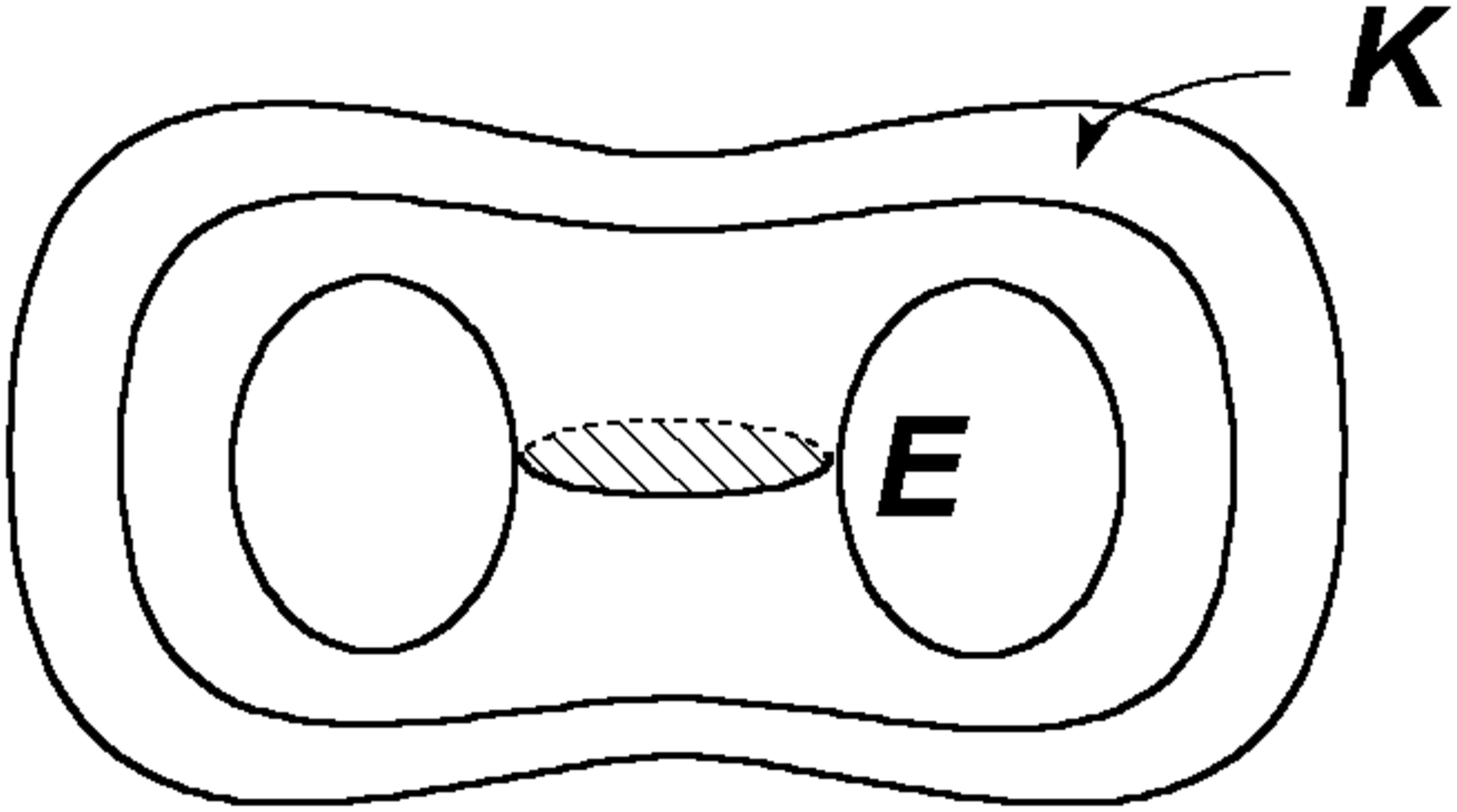}
\end{center}
\begin{center}
Figure D-1
\end{center}
\end{figure}

\medskip
\noindent{\bf Lemma D-1}\qua 
Let $D$ be a $K$--essential disk in $H$ 
such that $D \cap K$ consists of two points. 
Then there exists a $K$--boundary compressing disk
$\Delta$ for $D$. 

\medskip
\noindent{\bf Proof}\qua
Let $\tilde{D}$ be the lift of $D$ in $\tilde{H}$. 
Then, by Lemmas~A-1 and A-3, 
we see that $\tilde{D}$ is an essential annulus 
in a genus three handlebody $\tilde{H}$. 
Then $\tilde{D}$ is boundary compressible in $\tilde{H}$. 
Hence, by Lemma~A-2, we see that 
$D$ is $K$--boundary compressible in $H$. 
\qed 

\begin{figure}[ht!]\small
\begin{center}
\leavevmode
\epsfxsize=60mm
\epsfbox{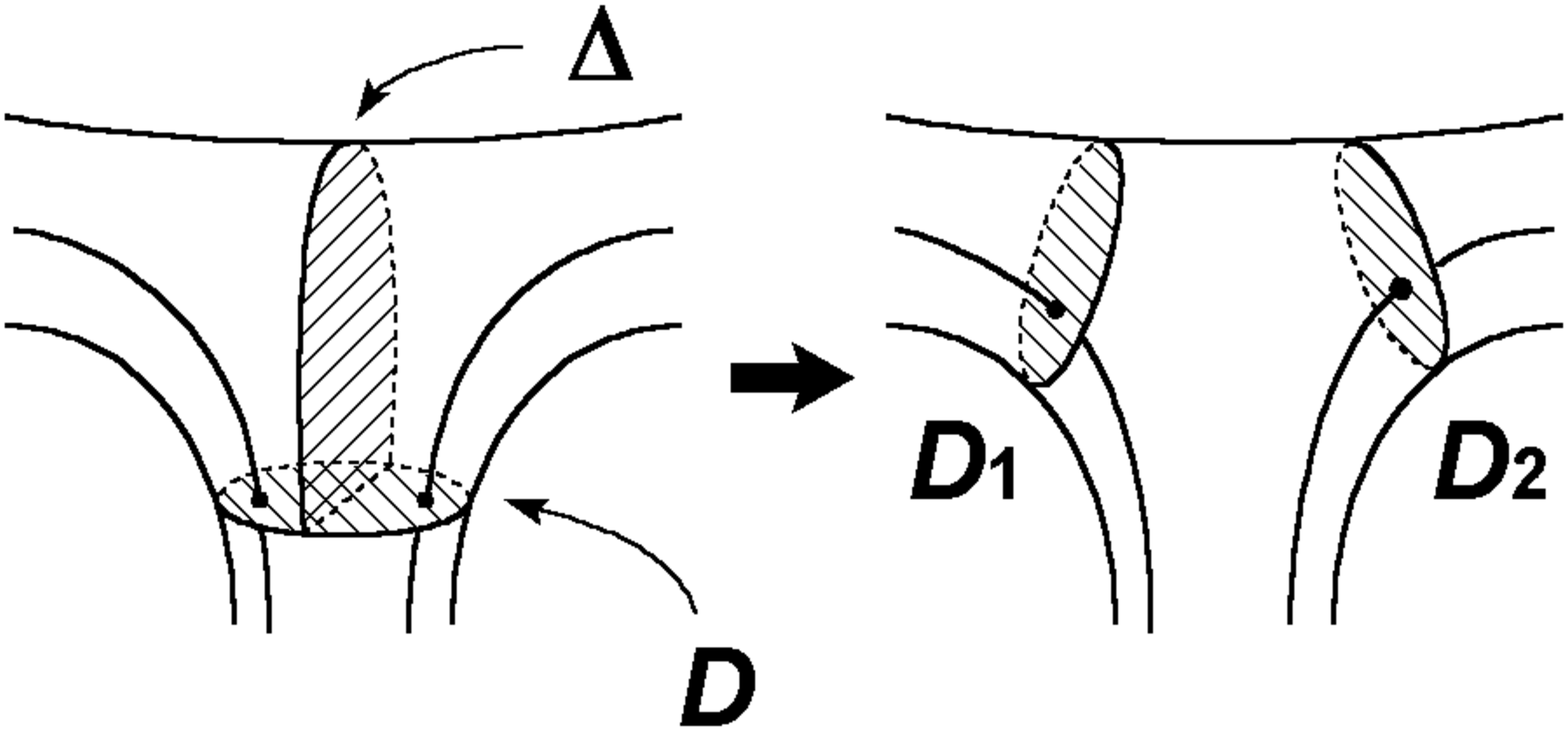}
\end{center}
\begin{center}
Figure D-2
\end{center}
\end{figure}

\medskip 
By Lemma~D-1, we obtain, 
by boundary compressing $D$ along $\Delta$, 
two $K$--compressing disks, 
say $D_1$ and $D_2$, for $\partial T$ 
such that $D_i \cap K$ consists of a point $(i=1,2)$. 

\medskip
\noindent{\bf Lemma D-2}\qua 
Let $D$, $D_1$, $D_2$ be as above. 
Suppose, moreover, that 
$D$ is separating in $H$. 
Then $D_1$ and $D_2$ are $K$--parallel in $H$, and, 
by cutting $(H, K)$ along $D_i$ ($i=1$ or $2$, say 1), 
we obtain a 1--string trivial arc 
in a solid torus, say $(T, \alpha )$. 
Moreover, $D_2$ is $\alpha$--boundary parallel in $T$. 

\begin{figure}[ht!]\small
\begin{center}
\leavevmode
\epsfxsize=40mm
\epsfbox{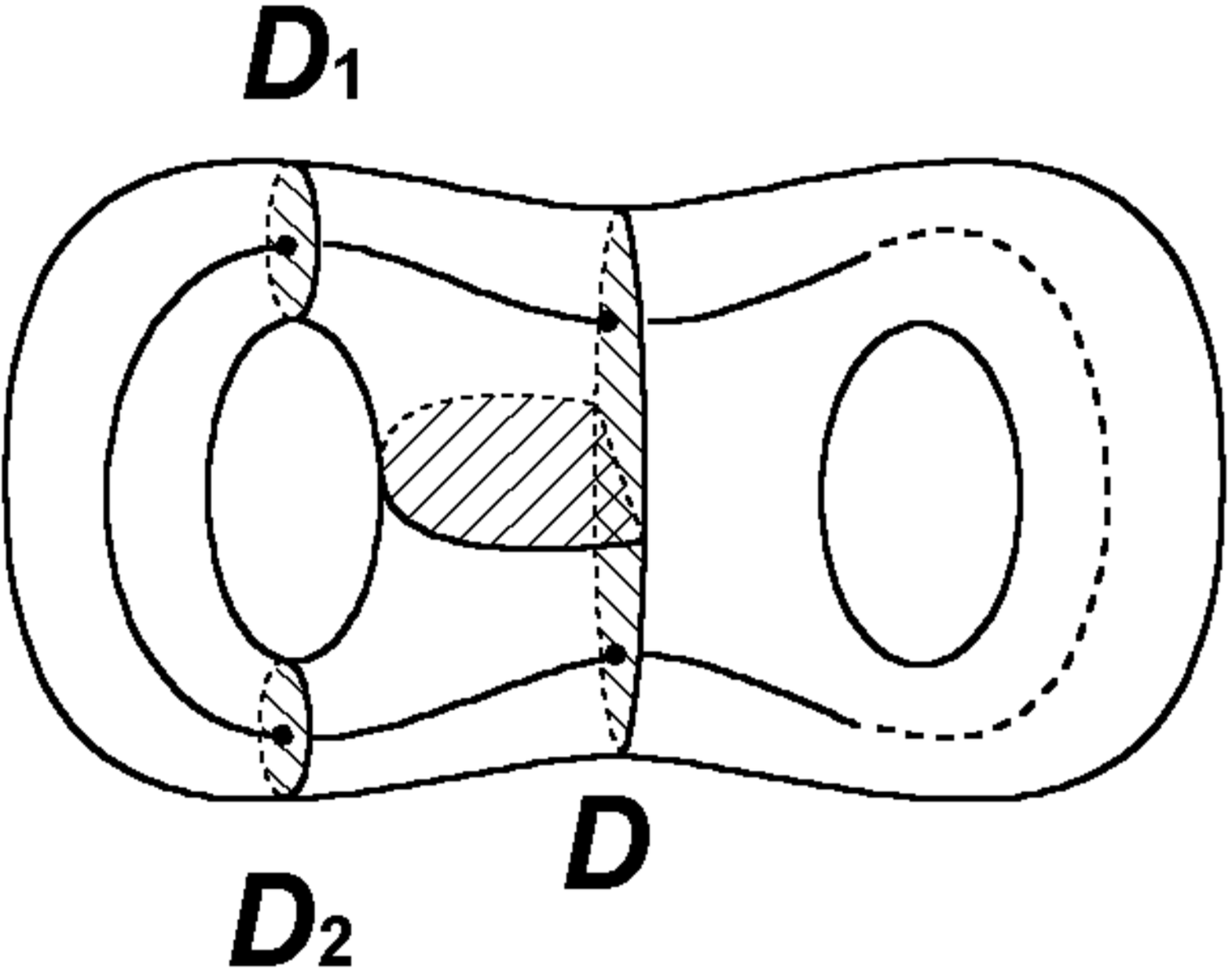}
\end{center}
\begin{center}
Figure D-3
\end{center}
\end{figure}

%\medskip
\noindent{\bf Proof}\qua
We note that $D$ separates $H$ into two solid tori $T_1$, $T_2$, 
where $D_1$, $D_2$ are properly embedded in $T_1$. 
Since each $D_i$ intersects $K$ in one point, 
$D_i$ is an essential disk of $T_1$, and this shows that 
$D_1$ and $D_2$ are parallel in $T_1$, and in $H$. 
Then, ${\Bbb Z}_2$--Smith conjecture shows that 
they are actually $K$--parallel. 
Then, by using ${\Bbb Z}_2$--equivariant loop theorem, 
we see that 
we obtain a 1--string trivial tangle in a solid torus
$(T, \alpha )$ , 
by cutting $(H, K)$ along $D_1$. 
Since $D_1$ and $D_2$ are $K$--parallel in $H$, 
we see that $D_2$ is $\alpha$--boundary parallel in $T$. 
\qed 

\medskip
\noindent{\bf Lemma D-3}\qua 
Let $D$ be as in Lemma~D-1. 
Suppose, moreover, that 
$D$ is non-separating in $H$. 
Then $D_1 \cup D_2$ is non-separating in $H$, and, 
by cutting $(H, K)$ along $D_1 \cup D_2$, 
we obtain a 2--string trivial tangle, say $(B, \beta )$. 
Moreover, $D$ is $\beta$--boundary parallel in $B$. 

\begin{figure}[ht!]\small
\begin{center}
\leavevmode
\epsfxsize=40mm
\epsfbox{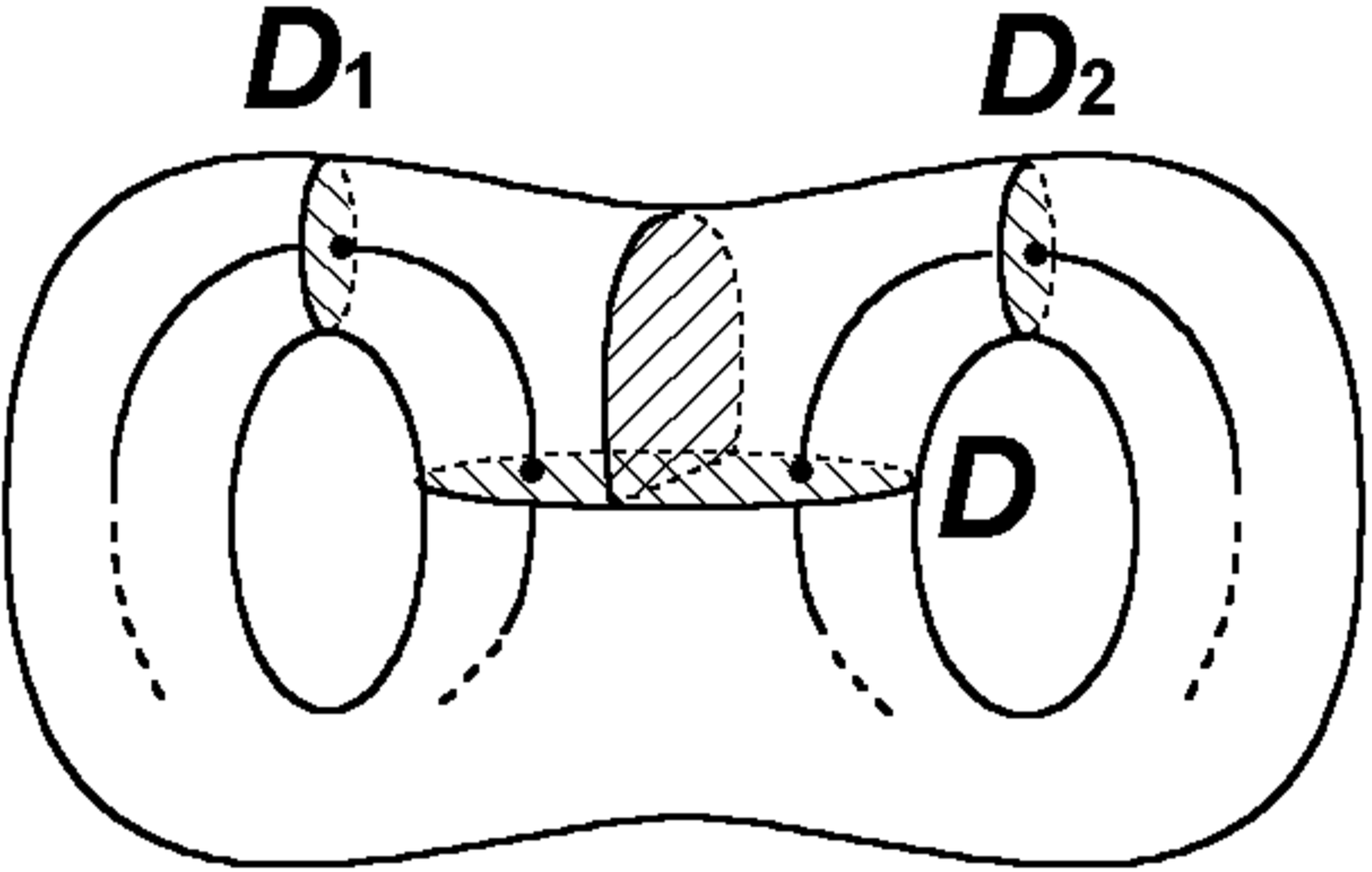}
\end{center}
\begin{center}
Figure D-4
\end{center}
\end{figure}

\medskip
\noindent{\bf Proof}\qua
%By the construction of $D_1$, $D_2$, 
%it is easy to see that 
Let $T$ be the solid torus obtained from $H$ by cutting along $D$. 
We may suppose that $D_1$ and $D_2$ are properly embedded 
in $T$. 
Since each $D_i$ intersects $K$ in one point, 
we see that $D_i$ is an essential disk in $T$. 
By the construction of $D_1$, $D_2$, 
we see that $D_1 \cup D_2$ separates the copies of 
$D$ in $T$. 
This shows that 
$D_1 \cup D_2$ is non-separating in $H$. 
Then, by cutting $(H, K)$ along $D_1 \cup D_2$, 
we obtain a 2--string tangle in a 3--ball, 
say $(B, \beta )$. 
Since $\tilde{H}$ is a genus three handlebody, 
we see that the 2--fold covering space 
of $B$ branched along $\beta$ is a 
solid torus. 
Hence, $(B, \beta )$ is a 2--string trivial tangle. 
By Lemma~B-1 (3), 
we see that $D$ is $\beta$--boundary parallel in $B$. 
\qed 

\medskip
\noindent{\bf Lemma D-4}\qua 
Let $D$, $D'$ be  
pairwise disjoint, pairwise parallel, 
non $K$--parallel, 
$K$--essential disks in $H$ 
such that $D \cap K$, and $D' \cap K$ 
consist of two points. 
Then there are two $K$--compressing disks 
$D^1$, $D^2$ for $\partial H$ such that 
$D^1 \cup D^2$ is non-separating in $H$ 
and is disjoint from $D \cup D'$, and, 
by cutting $(H, K)$ along $D^1 \cup D^2$, 
we obtain a 2--string trivial tangle, say $(B, \beta )$. 
Moreover $D$, $D'$ are $\beta$--boundary parallel in 
$(B, \beta )$, and, hence, 
$D \cup D'$ cobounds a 2--string trivial 
tangle in $(H, K)$. 

\medskip
\noindent{\bf Proof}\qua
Let $\Delta$ be a $K$--boundary compressing disk 
for $D \cup D'$. 
Without loss of generality, 
we may suppose that $\Delta \cap D \ne \emptyset$. 
We divide the proof into the following two cases. 

\medskip
\noindent
{\bf Case 1}\qua
$D$ and $D'$ are non-separating in $H$. 

\medskip
Let $D^1$, $D^2$ be the disks obtained from 
$D$ and $\Delta$ as in Lemma~D-3. 
Then, by the proof of Lemma~D-3, 
it is easy to see that $D^1 \cup D^2$ satisfies the 
conclusion of Lemma~D-4. 

\medskip
\noindent
{\bf Case 2}\qua
$D$ and $D'$ are separating in $H$. 

\medskip
Let $D^1$ be the disk corresponding to 
$D_1$ or $D_2$ in Lemma~D-2, 
and $(T, \alpha )$ the 1--string trivial arc 
in a solid torus $T$ obtained from $(H, K)$ by cutting 
along $D^1$. 
Then, by Lemma~B-2, 
we see that there exists an 
$\alpha$--compressing disk $D^2$ for $\partial T$ 
such that 
$D^2$ cuts $(T, \alpha )$ into a 2--string trivial tangle. 
Here we may suppose that 
$D^2$ is disjoint from the images of $D^1$ 
in $\partial T$, 
and, hence, 
we may regard that $D^2$ is properly embedded in $H$. 
Then $D^1 \cup D^2$ satisfies the 
conclusion of Lemma~D-4.

\begin{figure}[ht!]\small
\begin{center}
\leavevmode
\epsfxsize=70mm
\epsfbox{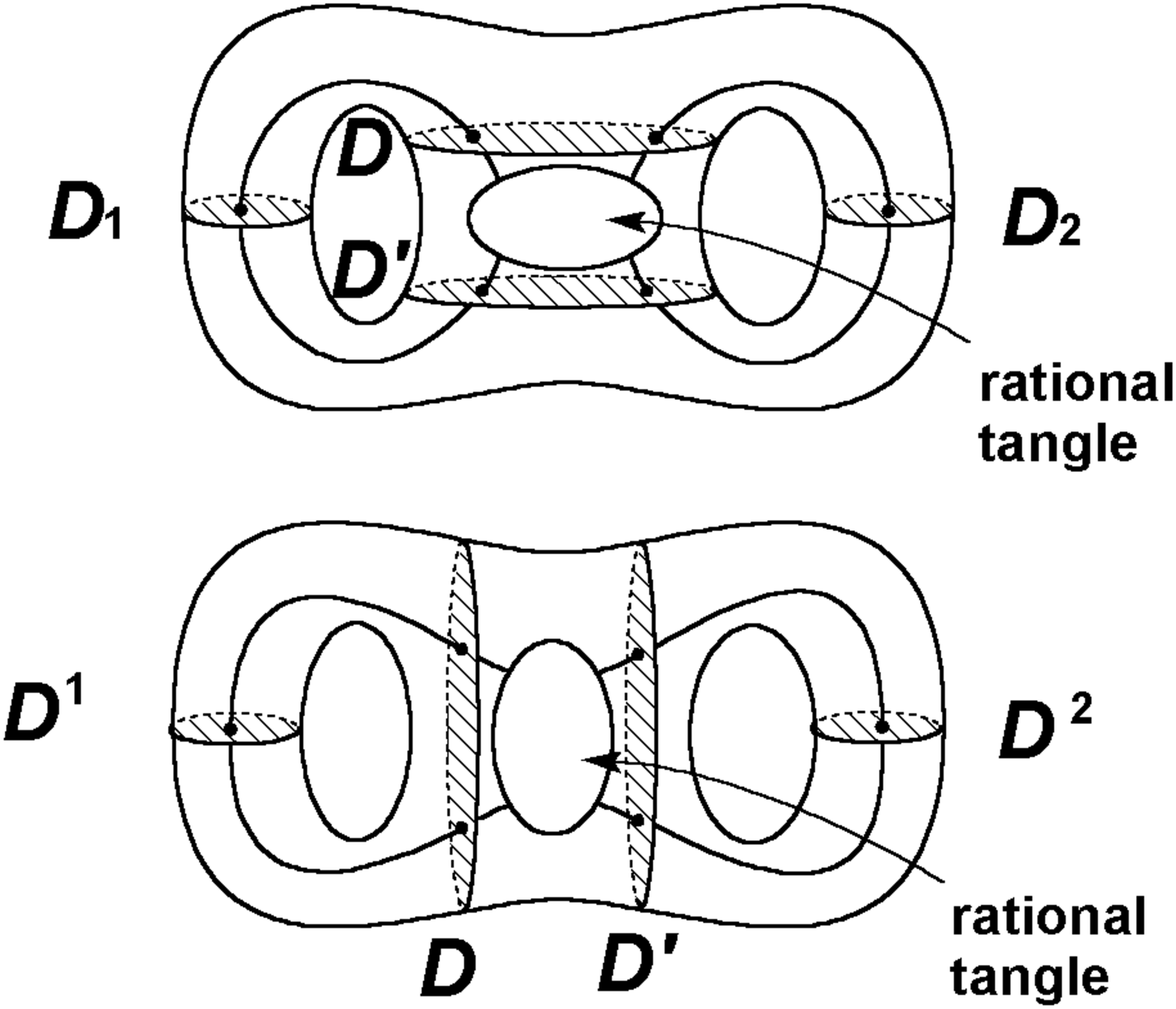}
\end{center}
\begin{center}
Figure D-5
\end{center}
\end{figure}

\end{document}